\documentclass[12pt]{article}
\usepackage{amsfonts}

\def\hybrid{\topmargin 0pt      \oddsidemargin 0pt
        \headheight 0pt \headsep 0pt
        \voffset=-0.5cm
        \textwidth 6.5in        
        \textheight 9in         
        \marginparwidth 0.0in
        \parskip 5pt plus 1pt   \jot = 1.5ex}
\catcode`\@=11
\def\marginnote#1{}

\newcount\hour
\newcount\minute
\newtoks\amorpm
\hour=\time\divide\hour by60
\minute=\time{\multiply\hour by60 \global\advance\minute by-\hour}
\edef\standardtime{{\ifnum\hour<12 \global\amorpm={am}%
        \else\global\amorpm={pm}\advance\hour by-12 \fi
        \ifnum\hour=0 \hour=12 \fi
        \number\hour:\ifnum\minute<10 0\fi\number\minute\the\amorpm}}
\edef\militarytime{\number\hour:\ifnum\minute<10 0\fi\number\minute}

\def\draftlabel#1{{\@bsphack\if@filesw {\let\thepage\relax
   \xdef\@gtempa{\write\@auxout{\string
      \newlabel{#1}{{\@currentlabel}{\thepage}}}}}\@gtempa
   \if@nobreak \ifvmode\nobreak\fi\fi\fi\@esphack}
        \gdef\@eqnlabel{#1}}
\def\@eqnlabel{}
\def\@vacuum{}
\def\draftmarginnote#1{\marginpar{\raggedright\scriptsize\tt#1}}
\def\draftlabel#1{{\@bsphack\if@filesw {\let\thepage\relax
   \xdef\@gtempa{\write\@auxout{\string
      \newlabel{#1}{{\@currentlabel}{\thepage}}}}}\@gtempa
   \if@nobreak \ifvmode\nobreak\fi\fi\fi\@esphack}
        \gdef\@eqnlabel{#1}}
\def\@eqnlabel{}
\def\@vacuum{}
\def\draftmarginnote#1{\marginpar{\raggedright\scriptsize\tt#1}}

\def\draft{\oddsidemargin -.5truein
        \def\@oddfoot{\sl preliminary draft \hfil
        \rm\thepage\hfil\sl\today\quad\militarytime}
        \let\@evenfoot\@oddfoot \overfullrule 3pt
        \let\label=\draftlabel
        \let\marginnote=\draftmarginnote
   \def\@eqnnum{(\theequation)\rlap{\kern\marginparsep\tt\@eqnlabel}%
\global\let\@eqnlabel\@vacuum}  }

\def\numberbysection{\@addtoreset{equation}{section}
        \def\theequation{\thesection.\arabic{equation}}}

\def\underline#1{\relax\ifmmode\@@underline#1\else
        $\@@underline{\hbox{#1}}$\relax\fi}

\def\titlepage{\@restonecolfalse\if@twocolumn\@restonecoltrue\onecolumn
     \else \newpage \fi \thispagestyle{empty}\c@page\z@
        \def\thefootnote{\fnsymbol{footnote}} }

\def\endtitlepage{\if@restonecol\twocolumn \else  \fi
        \def\thefootnote{\arabic{footnote}}
        \setcounter{footnote}{0}}  
\relax

\numberbysection
\hybrid

\def\beq{\begin{equation}}
\def\eeq{\end{equation}}
\def\bea{\begin{eqnarray}}
\def\eea{\end{eqnarray}}
\def\p{\partial}
\def\G{\Gamma}
\def\g{\gamma}
\def\s{\sigma}
\def\z{\zeta}
\def\L{{\cal L}}

\def\a{\alpha}
\def\b{\beta}
\def\e{\varepsilon}
\def\l{\lambda}
\def\f{\varphi}

\def\A{{\cal A}}

\def\V{{\cal V}}
\def\D{{\cal D}}
\def\F{{\cal F}}

\def\L{{\cal L}}

\def\O{{\cal O}}
\def\P{{\cal P}}

\def\dim{{\rm dim}}
\def\res{{\rm res}}

\def\wh{\widehat}

\def\Cee{\mathbb{C}}

\def\scrO{\mathcal{O}}

\def \matrix #1 {\left(\begin{array}{cc} #1 \end{array}\right)}
\newtheorem{theo}{Theorem}[section]

\newtheorem{cor}{Corollary}[section]
\newtheorem{lem}{Lemma}[section]

\begin{document}
\begin{titlepage}
\title{A characterization of Prym varieties}

\author{I.Krichever \thanks{Columbia University, New York, USA and
Landau Institute for Theoretical Physics, Moscow, Russia; e-mail:
krichev@math.columbia.edu. Research is supported in part by National Science
Foundation under the grant DMS-04-05519.}}

\date{April 15, 2006\footnote{the revised version}}

\maketitle

\begin{abstract} We prove that Prym varieties of algebraic curves with two smooth
fixed points of involution are exactly the indecomposable principally polarized
abelian varieties whose theta-functions provide explicit formulae for integrable
$2D$ Schr\"odinger equation.
\end{abstract}

\end{titlepage}

\section{Introduction}

The problem of the characterization of the Prymians among principally polarized abelian
varieties is almost as old as the famous Riemann-Schottky problem on the characterization
of the Jacobian locus. Until now despite all the efforts it has remained unsolved.
Analogs of quite a few geometrical characterizations of the Jacobians for the case
of Prym varieties are either unproved or known to be invalid
(see reviews \cite{taim1,shok} and references therein).

The first effective solution of the Riemann-Schottky problem was obtained by
T.Shiota (\cite{shiota}), who proved Novikov's conjecture: the Jacobains of curves
are exactly
the indecomposable principally polarized abelian varieties whose theta-functions
provide explicit solutions of the KP equation. Attempts to prove the analog of Novikov's
conjecture for the case of Prym varieties  were made in \cite{taim2, shiota2,bauv}.
In \cite{taim2} it was shown that Novikov-Veselov (NV) equation
provides local solution of the characterization problem. In \cite{shiota2,bauv}
the characterizations of the Prym varieties in terms of BKP and NV equations were proved
only under certain additional assumptions. Note, that in \cite{bauv} a counter example
showing that BKP equation has theta-functional solutions which do not correspond
to the Prym varieties was constructed.

The goal of this work is to solve the characterization problem of the Prym varieties using
the new approach proposed in the author's previous work \cite{kr-schot}, where it was shown
that  KP equation contains excessive information and the Jacobian locus can be
characterized in terms only of one of its {\it auxiliary} linear equations.

\begin{theo}(\cite{kr-schot}) An indecomposable symmetric matrix $B$ with positive
definite
imaginary part is the matrix of the $b$-periods of normalized holomorphic differentials on
a curve of genus g if and only if there exist $g$-dimensional vectors
$U\neq 0, V,A $ such that the equation
\beq\label{lax1}
\left(\p_y-\p_x^2+u\right)\psi=0
\eeq
is satisfied with
\beq\label{u2}
u=-2\p_x^2 \ln \theta (Ux+Vy+Z) \ \ {\rm and}\ \
\psi={\theta(A+Ux+Vy+Z)\over\theta(Ux+Vy+Z)}\, e^{p\,x+E\,y},
\eeq
where $p,E$ are constants.
\end{theo}
Here $\theta(z)=\theta(z|B), \ z=(z_1,z_2,\ldots, z_g)$ is the Riemann theta-function,
defined by the formula
\beq\label{teta1}
\theta(z)=\sum_{m\in \mathbb{Z}^g} e^{2\pi i(z,m)+\pi i(Bm,m)},\ \
(z,m)=m_1z_1+\cdots+m_gz_g
\eeq
The addition formula for the Riemann theta-function directly implies that
equation (\ref{lax1}) with $u$ and $\psi$ as in (\ref{u2}) is in fact equivalent to
the system of equations
\beq\label{gr}
\p_V\Theta[\e,0](A/2)-\p_U^2\Theta[\e,0](A/2)-2p\,\p_U\Theta[\e,0](A/2)+(E-p^2)
\Theta[\e,0](A/2)=0,
\eeq
where $\Theta[\e,0](z)=\theta[\e,0](2z|2B)$ are level two theta-functions with
half-integer characteristics $\e\in {1\over 2}Z_2^g$.

The characterization of the Jacobian locus given by Theorem 1.1 is stronger than
that given in terms of the KP equation (see details in \cite{flex}, where Theorem 1.1
was proved under the assumption that the closure $\langle A\rangle$ of the
subgroup  of $X$ generated by $A$ is irreducible). In terms of the
Kummer map,
\beq\label{kum}
\kappa:Z\in X\longmapsto
\{\Theta[\e_1,0](Z):\cdots:\Theta[\e_{2^g},0](Z)\}\in \mathbb{CP}^{2^g-1}\, ,
\eeq
the statement of Theorem 1.1 is equivalent to the characterization of
the Jacobians via flexes of the Kummer varieties, which is a particular case of
the trisecant conjecture, first formulated in \cite{wel1}.

The Prym variety of a smooth algebraic curve $\G$ with involution $\s:\G\longmapsto \G$
is defined as the {\it odd} subspace $\P(\G)\subset J(\G)$ of the Jacobian
with respect to the involution $\s^*: J(\G)\longmapsto J(\G)$ induced by
$\s$. It is principally polarized only if $\s$ has no fixed points or has two fixed
points $P_{\pm}$. In this work we consider only the second case.

Let $\G$ be a smooth algebraic curve with involution $\s$ having two fixed points $P_\pm$.
From the Riemann-Hurwitz formula it follows that if the genus of the
factor-curve $\G_0=\G/\s$ equals $g$, then the genus $\G$ is equal to $2g$. It is known that
on $\G$ there exists a basis of cycles $a_i, b_i$ with the canonical matrix of
intersections $a_i\cdot a_j=b_i\cdot b_j=0,\ a_i\cdot b_j=\delta_{ij},\ \ 1\leq i,j\leq 2g,$
such that $\s(a_k)=-a_{g+k},\ \s(b_k)=-b_{g+k}, 1\leq k\leq g$.
If $d\omega_i$ are normalized holomorphic differentials on $\G$, then the differential
$du_k=d\omega_k+d\omega_{g+k}$ are odd $\s^*(du_k)=-du_k$. By definition they are called
the normalized holomorphic Prym differentials.
The matrix of their $b$-periods
\beq\label{pi}
\Pi_{kj}=\oint_{b_k}du_j,\ \ 1\leq k,j\leq g\,,
\eeq
is symmetric, has positive definite imaginary part, and defines
the Prym theta-function $\theta_{Pr}=\theta(z|\Pi)$.

Before presenting our main result it is necessary to mention that
the Prym variety remains non-degenerate (compact) under certain degenerations of the curve.
No characterization of Prym varieties given in terms of equations for the matrix $\Pi$
of periods of the Prym differentials can single out the possibility of such degenerations.
An algebraic curve $\G$ that is smooth outside fixed points
$P_+,P_-,Q_1,Q_2,\ldots,Q_k$ of its involution $\s$, where  $P_{\pm}$ are smooth
and $Q_k$ are simple double points at which  $\s$ does not permute branches of
$\G$, will be denoted below by
$\{\G,\s,P_{\pm},Q_k\}$.

\begin{theo} An indecomposable principally polarized abelian variety $(X,\theta)$
is the Prym variety of a curve of type $\{\G,\s,P_{\pm},Q_k\}$,
if and only if there exist $g$-dimensional vectors $U\neq 0,V\neq 0,A$
such that one of the following equivalent conditions
holds:

$(A)$  The equation
\beq\label{sh}
\left(\p_x\p_t+u\right)\psi=0
\eeq
is satisfied with
\beq\label{anz}
u=2\p^2_{xt}\ln \theta (Ux+Vt+Z)+C, \ \ {\rm and}\ \
\psi={\theta(A+Ux+Vt+Z)\over \theta(Ux+Vt+Z)}\, e^{p\,x+Et},
\eeq
where $C,p,E$ are constants.

\medskip
$(B)$ The equations
\beq\label{k2}
\p^2_{UV}\Theta [\e,0](A/2)+p\,\p_V\Theta [\e,0](A/2)+E\,\p_U\Theta [\e,0](A/2)+
C\Theta [\e,0](A/2)=0
\eeq
are satisfied for all $\e\in {1\over 2}Z_2^g$.

\medskip
$(C)$ The equation
\beq\label{t}
\p_U\theta\,\p_V\theta\left(\p^2_U\p^2_V\theta\right)+
\p^2_U\theta\,\p^2_V\theta\left(\p_U\p_V\theta\right)-
\p^2_U\theta\,\p_V\theta\left(\p_U\p^2_V\theta\right)-
\p_U\theta\,\p^2_V\theta\left(\p^2_U\p_V\theta\right)|_{\Theta}=0
\eeq
is valid on the theta-divisor $\{Z\in\Theta: \theta(Z)=0\}$.
\end{theo}
The equivalence of $(A)$ and $(B)$ is a direct corollary of the addition formula for
the theta-function. The "if" part of $(A)$ follows from the construction
of integrable $2D$ Schr\"odinger operators given in \cite{nv}.
This construction is presented in the next section.

The statement $(C)$ is actually what we use for the proof of the theorem.
It is stronger than $(A)$. The implication $(A)\to (C)$ does not require
the explicit theta-functional form of $\psi$. It is enough to require only that
equation (\ref{sh}) with $u$ as in (\ref{anz}) has {\it local meromorphic} in $x$ (or $t$)
solutions which are holomorphic outside  the divisor $\theta(Ux+Vt+Z)=0$.

To put it more precisely, let us consider a function $\tau(x,t)$ which is a
holomorphic function of $x$ in
some domain where the equation $\tau(x,t)=0$ has a simple root $\eta(t)$. It turns out
that equation (\ref{sh}) with the potential $u=2\p^2_{xt}\ln \tau+C$, where $C$ is a constant,
has a meromorphic solution in $D$, if this root satisfies the equation
\beq\label{sys}
\ddot\eta v-\dot \eta\dot v+2\dot \eta^2w=0,
\eeq
where $v=v(t), \ w=w(t)$ are the first coefficients of the Laurent expansion of $u$ at $\eta$
\beq\label{e1}
u(x,t)={2\dot \eta\over (x-\eta)^2}+v+w (x-\eta)+\cdots
\eeq
and "dots" stand for $t$-derivatives. Straightforward but tedious computations with
expansion of $\theta$ at the generic points of
its divisor $\Theta$ show that equation
(\ref{sys}) in the case when $\tau=\theta(Ux+Vt+Z)$ is equivalent to equation (\ref{t}).

Note that equations (\ref{sys}) are analogues of the equations derived in \cite{flex}
and called in \cite{kr-schot} the formal Calogero-Moser system. In a similar way, if
we represent an entire function $\tau$ as a product
\beq\label{roots}
\tau(x,t)=c(t)\prod_i(x-x_i(t)),
\eeq
then equation (\ref{sys}) takes the form
\beq\label{cm}
\sum_{j\neq i}\left[{\ddot x_i\dot x_j-\dot x_i\ddot x_j\over(x_i-x_j)^2}-
{2\dot x_i\dot x_j(\dot x_i+\dot x_j)\over(x_i-x_j)^3}\right]=0.
\eeq
At the moment the only reason for presenting equations (\ref{cm})
is to show that in the case when $\tau$ is a rational, trigonometric or
elliptic polynomial  the system (\ref{sys}) gives well-defined equations of
motion for a multi-particle system. \footnote{A.Zotov noticed that equation
(\ref{cm}) is equivalent to the equations
$\ddot x_i=2\sum_{j\neq i} \dot x_i\dot x_j/(x_i-x_j)$ and, therefore, can be regarded
as a limiting case of the Ruijesenaars-Schneider system.}

At the beginning of  section 3 we derive equation (\ref{sys}) and show that equation
(\ref{t}) is sufficient for the {\it local} existence of wave solutions of (\ref{sh})
having the form
\beq\label{ps}
\psi_{\pm}(x,t,k)=
e^{kt_{\pm}}\left(1+\sum_{s=1}^{\infty}\xi^{\pm}_s(x,t)\,k^{-s}\right)\,, \ t_+=x,\,t_-=t,
\eeq
and such that
\beq\label{xi22}
\xi^{\pm}_s={\tau^{\pm}_s(Ux+Vt+Z,t_{\mp})\over\theta (Ux+Vt+Z)}\ , \ Z\notin \Sigma_{\pm},
\eeq
where $\tau^{\pm}_s(Z,t_{\mp})$, as a function of $Z$, is holomorphic in some open domain
in $\mathbb C^g$. Here and below $\Sigma_{\pm}\subset \Theta$ are subsets of the theta-divisor
invariant under the shifts along constant vector fields $\p_U$ or $\p_V$, respectively.

The coefficients  $\xi_s^{\pm}$ of the wave solutions are defined recurrently by the
equations $\partial_{\mp}\xi^{\pm}_{s+1}=-\p_{xt} \xi_s-u\xi_s$. The local existence of
meromorphic solutions requires vanishing of the residues of the nonhomogeneous terms.
That is controlled by equation (\ref{sys}). At the local level
the main problem is to find the translational invariant normalization of $\xi^{\pm}_s$ which defines
wave solutions uniquely up to a $(x,t)$-independent factor.

Following the ideas of \cite{kp2} and \cite{kr-schot} we fix such a normalization using
extensions of $\xi_s^{\pm}$ along the affine subspaces $Z+{\mathbb C}^d_{\pm}$, where
${\mathbb C}^d_{\pm}$ are universal covers of the abelian subvarieties $Y_{\pm}\subset X$
which are closures of the subgroups $Ux$ and $Vt$ in $X$, respectively.
The corresponding wave solutions are called $\l$-periodic.

In the last section we show that for each $Z\notin \Sigma_{\pm}$ a local
$\l$-periodic wave solution is the common eigenfunction of a commutative ring $\A_{\pm}^Z$
of ordinary differential operators. The coefficients of these operators are independent of
ambiguities in the construction of $\psi$. The theory of commuting differential operators
\cite{kr1,kr2,ch1,ch2} implies then that the correspondence
$Z\longmapsto \A_{\pm}^Z$ defines a map
$j$ of $X\setminus\Sigma_{\pm}$ into the space $\overline{\rm Pic(\G)}$
of torsion-free rang 1 sheafs $\F$ on $Z$-independent spectral curve of $\A^Z$.
That allows us to make the next crucial step and prove the global existence of
the wave function. The global existence of the
wave function implies that for the generic $Z\notin \Sigma_{\pm}$ the orbit of $\A^Z$
under the NV flows defines an imbedding $i_Z$ of the Prym
variety $\P(\G)$ of the spectral curve into $X$.
Therefore, the Prym variety is compact. That implies the explicit description
of possible types of singular points of $\G$. The final step in the proof of the main theorem
is to show that there are not singular points
of the multiplicity bigger then 2.

\section{Integrable 2D Schr\"odinger operators}

In this section we present necessary facts from the theory of
integrable $2D$ Schr\"odinger equations and related hierarchies.

Let $\G$ be a smooth algebraic curve of genus $g$ with fixed local coordinates $k^{-1}_{\pm}$
at punctures $P_{\pm}$ and let $t^{(\pm)}=\{t_i^{(\pm)}\}$ be finite
sets of complex variables. Then according to the general construction of the multi-point
Baker-Akhiezer functions (\cite{kr1,kr2}) for each non-special effective
divisor $D=\{\g_1,\ldots,\g_g\}$ of degree $g$ there exists a unique function
$\psi_0(t^{(+)},t^{(-)}, P)$, which, as a function of the variable $P\in\G$, is meromorphic on
$\G\setminus P_{\pm}$, where it has poles at $\g_s$ of degree not greater then the multiplicity
of $\g_s$ in $D$. In the neighborhood of $P_{\pm}$ the function $\psi_0$ has the form
\beq\label{ps1}
\psi_0=
e^{\sum_ik^it^{(\pm)}_i}\left(\sum_{s=0}^{\infty}\xi^{\pm}_s(t)\, k^{-s}\right),
\ \ \xi_0^+=1,
\eeq
where $k=k_{\pm}^{-1}(P)$ and $t=\{t^{(+)},t^{(-)}\}$.

The uniqueness of $\psi_0$ implies that for each positive integer $n$ there exists
unique differential operators $B_n^{\pm}$ in the variables $t_1^{\pm}$
\beq\label{ln}
B_n^{\pm}=\p_{\pm}^n+\sum_{i=0}^{n-1}v_{n,i}^{\pm}(t)\p_{\pm}^{n-i},\ \
\p_{\pm}=\p/\p t_1^{\pm},
\eeq
such that
\beq\label{t1}
\left(\p_{t_n^{\pm}}-B_n^{\pm}\right)\psi_0=0.
\eeq
Equations (\ref{t1}) directly imply
\beq\label{t2}
\left[\p_{t_n^{\pm}}-B_n^{\pm},\p_{t_m^{\pm}}-B_m^{\pm}\right]=0.
\eeq
In other words, the operators $B_n^{\pm}$ satisfy zero-curvature equations which define
{\it two} copies of the KP hierarchy with respect to the times $t_n^{\pm}$.

The {\it two-point} Baker-Akhiezer function with {\it separated variables} was introduced in
\cite{dkn} where it was proved that in addition to (\ref{t1}) it satisfies the
equation
\beq\label{t3}
H\psi_0=\left(\p_+\p_-+w\p_++u\right)\psi_0=0\,,
\eeq
where
\beq\label{t4}
w=-\p_-\ln\xi_0^-,\ \ u=-\p_+\p_-\ln \xi_1^+\,.
\eeq
The operator $H$ defined in the left hand side of (\ref{t3}) "couples" two copies
of the KP hierarchy corresponding to the punctures $P_{\pm}$ via the equation
\beq\label{t5}
\left[\p_{t_n^{+}}-B_n^{+},\p_{t_m^{-}}-B_m^{-}\right]=D_{nm}H.\
\eeq
The sense of (\ref{t5}) is as follows. Each differential operator $ \D $ in the two
variables $t^{\pm}_1$ can be uniquely represented in the form
\beq
\D = D H+ D^+ + D^-, \label{t6}
\eeq
where $D^{\pm}$ are ordinary differential operators in the variables
$t^{\pm}_1$, respectively. The equation (\ref{t5}) is just the statement that
the second and the third terms in the corresponding representation of the left hand side
of (\ref{t5}) are equal to zero. This
implies $n+m-1$ equations on $n+m-1$ unknown functions (the coefficients of
$ B^+_n$ and $B^-_m$ ). Therefore, the operator equation (\ref{t6}),
is equivalent to the well-defined system of non-linear partial
differential equations.

Explicit theta-functional formulae for the solutions of these equations follow from
the theta-functional formula for the Baker-Akhiezer function
\beq
\psi_0=
{\theta(A(P)+\sum_{i} \left(U_{i}^+t_{i}^++U_{i}^-t_{i}^-\right)+Z)\,\theta(A(P_+)+Z)
\over \theta(A(P_+)+\sum_{i} \left(U_{i}^+t_{i}^++U_{i}^-t_{i}^-\right)+Z)\,\theta(A(P)+Z)}
 \ e^{\ \sum_{i} \left(t_{i}^+\Omega_{i}^+(P)+t_{i}^-\Omega_{i}^-(P)\right)} \label{2.101}
\eeq
Here:

a) $\theta(z)=\theta(z|B)$ is the Riemann theta-function defined
by the matrix $B$ of $b$-periods of normalized holomorphic differentials $d\omega_k$ on $\Gamma$.

b) $\Omega_{i}^{\pm}(P)=\int^P d\Omega_{i}^{\pm}$ is the Abelian integral
corresponding to the normalized,
$\oint_{a_k} d\Omega_{i}^{\pm}=0,$
meromorphic differential on $\Gamma$ with the only pole of the form
\beq
d\Omega_{i}^{\pm}=dk_{\pm}^i(1+O(k_{\pm}^{-i-1})) \label{2.104a}
\eeq
at the puncture $P_{\pm}$;

c) $2\pi iU_{j}^{\pm}$ is a vector of $b$-periods of the differential
$d\Omega_{j}^{\pm}$ with the coordinates
\beq
U_{j,k}^{\pm}={1\over 2\pi i} \oint_{b_k} d\Omega_{j}^{\pm}; \label{2.105}
\eeq

d) $A(P)$ is the Abel  transform, i.e. it is a vector with the coordinates
$A_k(P)=\int^P d\omega_k$;

e) $Z$ is an arbitrary vector (it corresponds to the divisor of poles of
Baker-Akhiezer function).

Taking the evaluation at $P_-$ and the expansion at $P_+$ of the regular factor in
(\ref{2.101}), one gets theta-functional formulae for the coefficients (\ref{t4}) of the
corresponding
2D Schr\"odinger operator:
\begin{eqnarray}\label{t7}
w&=&-\p_-\ln\left({\theta(A(P_-)+\sum_{i} \left(U_{i}^+t_{i}^++U_{i}^-t_{i}^-\right)+Z)
\over \theta(A(P_+)+\sum_{i} \left(U_{i}^+t_{i}^++U_{i}^-t_{i}^-\right)+Z)}\right)\,,\\
u&=&\p_+\p_-\ln \left(\theta(A(P_+)+\sum_{i} \left(U_{i}^+t_{i}^++U_{i}^-t_{i}^-\right)+Z)
\right)+C\,,\label{t8}
\end{eqnarray}
where the constant $C$ is equal to $C=\res_{P_+}\Omega_-d\Omega_+$.
Note, that the second factors in the numerator and denominator of the formula (\ref{2.101})
are $t$-independent. Therefore, the function $\psi$ given by the following formula
\beq\label{unnor}
\psi=
{\theta(A(P)+\sum_{i} \left(U_{i}^+t_{i}^++U_{i}^-t_{i}^-\right)+Z)
\over \theta(A(P_+)+\sum_{i} \left(U_{i}^+t_{i}^++U_{i}^-t_{i}^-\right)+Z)}
 \ e^{\ \sum_{i} \left(t_{i}^+\Omega_{i}^+(P)+t_{i}^-\Omega_{i}^-(P)\right)}
\eeq
is a solution of the same linear equations as $\psi_0$. Below $\psi$ given by (\ref{unnor})
 will be called {\it non-normalized} Baker-Akhiezer function.

{\bf Potential operators}. From now on we will consider only {\it potential} Schr\"odinger operators
$H=\p_+\p_-+u$. The reduction of the above described algebraic-geometrical construction
to the potential case was found in \cite{nv}. The corresponding algebraic-geometrical
data are singled out by the following constraints:

$(i)$ {\it The curve $\G$ should be a curve with involution $\sigma:\G \to \G$
which has {\it two} fixed points $P_{\pm}$}.

$(ii)$ {\it The equivalence class $[D]\in J(\G)$
of the divisor $D$ should satisfy the equation}
\beq
[D]+[\s(D)]=K+P_++P_-\in J(\G),\label{p1}
\eeq
{\it where $K$ is the canonical class, i.e. the equivalence class of the
zero-divisor of a holomorphic differential on $\G$}.

Equation (\ref{p1}) is equivalent to the condition that the divisor
$D+\s(D)$ is the zero divisor of a meromorphic differential $d\Omega$ on $\G$ with
simple poles at the punctures $P_{\pm}$. The differential $d\Omega$ is even with
respect to the involution and descends to a meromorphic differential on the
factor-curve $\G_0=\G/ \sigma$. The projection $\pi :\Gamma \longmapsto \Gamma_0=\Gamma / \sigma$
represents $\Gamma$ as a two-sheet covering of $\Gamma_0$ with $2$
branch points $P_{\pm}$. In this realization the involution $\sigma$ is
a permutation of the sheets. For $P\in \G$ we denote the point
$\sigma(P)$ by $P^{\s}$. From the Riemann-Hurwitz formula it follows that the
genus $g$ of $\G$ equals $g=2g_0$, where $g_0$ is the genus of $\Gamma_0$.
Note that the divisors that satisfy (\ref{p1}) are parameterized by the points $Z_0$ of
the Prym variety $\P(\G)\subset J(\G)$.

\begin{theo} (\cite{nv}). Let a smooth algebraic curve $\G$ and
an effective divisor $D$ satisfy the constraints $(i), (ii)$. Let $k^{-1}_{\pm}(P)$
be odd local coordinates in the neighborhoods of the fixed points $P_{\pm}$, i.e.
$k_{\pm}(P)=-k_{\pm}(\sigma(P))$, and let all the even times vanish, i.e. $t_{2i}^{\pm}=0$.
Then the corresponding $2D$ Schr\"odinger operator is potential, i.e. $w=0$.
\end{theo}
In \cite{nv}) it was also found that for the potential operators
the formulae (\ref{2.101}) and (\ref{t8}) can be expressed in terms of
the Prym theta-function.
For further use it is enough to present these formulae for the case of only two
nontrivial variables $x=t_1^+,\ t=t_1^-$:
\beq\label{nv1}
\psi=
{\theta_{Pr}(A^{Pr}(P)+Ux+Vt+Z)
\over \theta_{Pr}(A^{Pr}(P_+)+Ux+Vt+Z)}
 \ e^{\ x\Omega_1^++t\Omega_1^-}
\eeq
\beq\label{nv2}
u=2\p_+\p_-\ln \,\theta_{Pr}(A^{Pr}(P_+)+xU+tV+Z)+C\,.
\eeq
Here $A^{Pr}:\G\longmapsto \P(\G)$ is the Abel-Prym map defined by the Prym differentials,
i.e. $A^{Pr}(P)$ is a vector with the coordinates $A^{Pr}_k(P)=\int^P du_k$.

In \cite{kr-spec,taim3} it was proved that for the case of smooth periodic
potentials $u(x,t)$ (considered as a function of real variables $x,t$) the conditions found
by Novikov and Veselov are {\it sufficient and necessary}.

\section{$\l$-periodic wave solutions}

To begin with, let us show that equations (\ref{sys}) are the
necessary condition of the existence of a meromorphic solution to equation (\ref{sh}).

Let $\tau(x,t)$ be a smooth $t$-parametric family of holomorphic functions of
the variable $x$ in some open domain $D\subset \mathbb C$. Suppose that in $D$ the function
$\tau$ has a simple zero,
\beq\label{xi}
\tau(\eta(t),t)=0, \tau_x(\eta(t),t)\neq 0.
\eeq
\begin{lem}
If equation (\ref{sh}) with the potential
$u=2\p_{xt}^2\ln \tau(x,t)+C$, where $C$ is a constant,
has a meromorphic solution $\psi_0(x,t)$, then
equation (\ref{sys}) holds.
\end{lem}
{\it Proof.}
Consider the Laurent expansions of $\psi_0$ and $u$ in the neighborhood of $\eta$:
\beq\label{ue}
u={2\dot \eta\over (x-\eta)^2}+v+w(x-\eta)+\ldots
\eeq
\beq\label{psie}
\psi_0={\a\over x-\eta}+\b+\g(x-\eta)+\ldots
\eeq
(All the coefficients in these expansions are smooth functions of the variable $t$).
Substitution of (\ref{ue},\ref{psie}) into (\ref{sh})
gives a system of equations. The first three of them
are
\beq\label{eq1}
\dot\a-2\dot\eta\b=0,
\eeq
\beq\label{eq2}
2\dot\eta\g+\a v=0,
\eeq
\beq\label{eq3}
\dot\g+v\b+\a w=0.
\eeq
Taking the $t$-derivative of the second equation and using two others we get
(\ref{sys}).

Let us show that equations (\ref{sys}) are sufficient for the existence of
meromorphic wave solutions.
\begin{lem}
Suppose that equation (\ref{sys}) for the zero of $\tau(x,t)$ holds. Then equation
(\ref{sh}) has wave solutions of the form
\beq\label{w1}
\psi=e^{kt}\left(1+\sum_{s=1}^{\infty} \xi_s(x,t)k^{-s}\right)
\eeq
such that the coefficients $\xi_s$ have simple poles at $\eta$ and
are holomorphic everywhere else in $D$.
\end{lem}
{\it Proof.} Substitution of (\ref{w1}) into (\ref{sh}) gives a recurrent system of
equations
\beq\label{xis}
\xi_{s+1}'=-\p^2_{xt}\xi_s-u\xi_s.
\eeq
We are going to prove by induction that this system has meromorphic solutions with
simple poles at $\eta$.

Let us expand $\xi_s$ at $\eta$:
\beq\label{5}
\xi_s={r_s\over x-\eta}+r_{s0}+r_{s1}(x-\eta)+\cdots
\eeq
Suppose that $\xi_s$ is defined, and equation (\ref{xis}) has a meromorphic solution.
Then the right hand side of (\ref{xis}) has the zero residue at $x=\eta$, i.e.,
\beq\label{res}
{\rm res}_{\eta}\left(\p^2_{xt}\xi_s+u\xi_s\right)=vr_s+2\dot\eta r_{s1}=0.
\eeq
We need to show that the residue of the next equation also vanishes.
From (\ref{xis}) it follows that the coefficients of the Laurent expansion for $\xi_{s+1}$
are equal to
\beq\label{6}
r_{s+1}=-\dot r_s+2\dot \eta r_{s0},
\eeq
\beq\label{7}
 r_{s+1,1}=-vr_{s0}-wr_s-\dot r_{s1}\,.
\eeq
These equations and equation (\ref{res}) imply
\beq
-(vr_{s+1}+2\dot\eta r_{s+1,1})=2\dot\eta w r_s+ v\dot r_s+2\dot\eta \dot r_{s,1}=
\p_t\left(vr_s+2\dot \eta r_{s,1}\right)-
\left(\dot v-{\ddot \eta\over \dot \eta}v-2\dot\eta w\right)r_s=0,
\eeq
and the lemma is proved.

Our next goal is to fix a {\it translation-invariant} normalization of $\xi_s$
which defines wave functions uniquely up to a $(x,t)$-independent factor.
It is instructive to consider first the case of the periodic potentials $u(x+1,t)=u(x,t)$
(compare with \cite{kp2}).

Equations (\ref{xis}) are solved recursively by the formulae
\beq
\xi_{s+1}(x,t)=c_{s+1}(t)+\xi_{s+1}^0(x,t)\,,\label{kp1}
\eeq
\beq\label{kp2}
\xi_{s+1}^0(x,t)=-\p_t \xi_s-\int_{x_0}^x u\xi_s\,dx\, ,
\eeq
where $c_s(t)$ are {\it arbitrary} functions of the variable $t$.
Let us show that the periodicity condition $\xi_s(x+1,t)=\xi_s(x,t)$
defines these functions uniquely up to constants.
Assume that $\xi_{s-1}$ is known  and satisfies the condition that the corresponding
function $\xi_s^0$ is periodic.
The choice of the function $c_s(t)$ does not affect the periodicity property of
$\xi_s$, but it does affect the periodicity in $x$ of the function
$\xi_{s+1}^0(x,t)$. In order to make  $\xi_{s+1}^0(x,t)$ periodic,
the function $c_s(t)$ should satisfy the linear differential equation
\beq\label{kp4}
\p_t c_s(t)+\int_{x_0}^{x_0+1}
u(x,t)\,(c_s(t)+\xi_s^0(x,t))\,dx\ =0.
\eeq
This defines $c_s$ uniquely up to a constant.

In the general case, when $u$ given by (\ref{anz}) is quasi-periodic, the normalization of
the wave functions is defined along the same lines.

Let $\Theta_1$ be defined by the equations $\Theta_1=\{Z:\theta(Z)=\partial_U\theta(Z)=0\}$,
where $\p_U$ is a constant vector-field on ${\mathbb C}^g$, corresponding to the vector $U$ in
(\ref{anz}). The $\p_U$-invariant subset $\Sigma$ of $\Theta_1$ will be called the
{\it singular locus}.

Consider the closure $Y_U=\langle Ux\rangle$ of the group
$Ux$ in $X$. Shifting $Y_U$ if needed, we may assume, without loss of generality, that
$Y_U$ is not in the singular locus, $Y_U\notin\Sigma$. Then, for a sufficiently small
$t$, we have $Y_U+Vt\notin\Sigma$ as well.
Consider the restriction of the theta-function onto the affine subspace
$\mathbb C^d+Vt$, where $\mathbb C^d=\pi^{-1}(Y_U)$, and
$\pi: \mathbb C^g\to X=\mathbb C^g/\Lambda$ is the universal cover of $X$:
\beq\label{ttt1}
\tau (z,t)=\theta(z+Vt), \ \ z\in \mathbb C^d.
\eeq
The function $u(z,t)=2\partial_U\p_t\ln \tau+C$ is periodic with respect to the lattice
$\Lambda_U=\Lambda\cap \mathbb C^d$ and, for fixed $t$, has a double pole along the divisor
$\Theta^{\,U}(t)=\left(\Theta-Vt\right)\cap \mathbb C^d$.

\begin{lem} Let equation (\ref{t}) hold and let
$\l$ be a vector of the sublattice $\Lambda_U=\Lambda\cap \mathbb C^d\subset \mathbb C^g$.
Then:

(i) equation (\ref{sh}) with the potential $u(Ux+z,t)$
has a wave solution of  the form $\psi=e^{kt+bxk^{-1}}\phi(Ux+z,t,k)$
such that the coefficients $\xi_s(z,t)$ of the formal series
\beq\label{psi2}
\phi(z,t,k)=1+\sum_{s=1}^{\infty}\tilde\xi_s(z,t)\, k^{-s}
\eeq
are $\l$-periodic meromorphic functions of the variable $z\in \mathbb C^d$
with a simple pole at
the divisor $\Theta^U(t)$, i.e.
\beq\label{v1}
\tilde\xi_s(z+\l,t)=\tilde\xi_s(z,t)={\tau_s(z,t)\over \tau(z,t)}\, ;
\eeq

(ii) $\phi(z,t,k)$ is unique up to a factor $\rho(z,k)$ that is $t$-independent,
$\partial_U$-invariant and holomorphic in $z$,
\beq\label{v2}
\phi_1(z,t,k)=\phi(z,t,k)\rho(z,k), \ \partial_U\rho=0.
\eeq
\end{lem}
{\it Proof.} The functions $\tilde\xi_s(z)$ are defined recursively by the equations
\beq\label{xis1}
\p_U\tilde\xi_{s+1}=-\p_U\p_t\tilde\xi_s-(u+b)\tilde\xi_s-b\p_t\tilde\xi_{s-1}.
\eeq
A particular solution of the first equation $\partial_U\tilde
\xi_1=-u-b$ is given by the formula
\beq\label{v5}
\tilde \xi_1^0=-2\p_t\ln \tau -(l,z)\ (b+C),
\eeq
where $(l,z)$ is a linear form on $\mathbb C^d$ given by the scalar product of $z$
with a vector $l\in \mathbb C^d$ such that $(l,U)=1$. By definition, the vector $\l$
is in $Y_U$. Therefore, $(l,\l)\neq 0$. The periodicity condition for
$\tilde\xi_1^0$ defines the constant $b$, which depends only on a choice of the
lattice vector $\l$.
From the monodromy properties of $\theta$ it follows that without loss of generality we
may assume that $\l$ is chosen such that the corresponding
constant $b$ is not equal to zero, i.e.
\beq\label{v6}
b=-C+(l,\l)^{-1}(2\p_t\ln \tau(z,t)-2\partial_t\ln \tau(z+\l,t))\neq 0\,,
\eeq
Note, that the second factor in (\ref{w1}) and the series $\phi$ in (\ref{psi2})
differ by the factor $e^{bxk^{-1}}$, which does not affect the results of the previous
lemma. Therefore, equations (\ref{sys}) are sufficient for the local solvability of
(\ref{xis1}) in any domain, where $\tau(z+Ux,t)$ has simple zeros, i.e. outside
the set $\Theta_1^{\,U}(t)=\left(\Theta_1-Vt\right)\cap \mathbb C^d$.
Recall that $\Theta_1= \Theta\cap \partial_U\Theta$.
This set  does not contain a $\partial_U$-invariant line because  such line is dense in
$Y_U$. Therefore, the sheaf $\V_0$ of $\partial_U$-invariant meromorphic functions on
$\mathbb{C}^d\setminus \Theta^{\,U}_1(t)$ with poles along the divisor $\Theta^{\,U}(t)$
coincides with the sheaf of holomorphic $\partial_U$-invariant functions.
That implies the vanishing of $H^1({C}^d\setminus \Theta^{\,U}_1(t), \V_0)$ and the existence
of global meromorphic solutions $\xi_s^0$ of (\ref{xis1}) which have a simple pole at
the divisor $\Theta^{\,U}(t)$ (see details in \cite{shiota, ac}).

Let us assume, as in the example above, that a $\l$-periodic solution
$\tilde\xi_{s-1}$ is known  and that it satisfies the condition that there exists
a $\l$-periodic solution $\tilde\xi_s^0$ of the next equation such that the equation
\beq\label{ap1}
\p_U\chi_{s}=-(u+b)\tilde\xi_s^0-b\p_t\tilde\xi_{s-1}
\eeq
has a $\l$-periodic solution. If $\tilde\xi_s^0$ and a particular solution
$\chi^*_{s}$ of (\ref{ap1}) are fixed,
then $\tilde\xi^*_{s+1}=\p_t\tilde\xi_s^0+\chi_{s}^*$ is a
$\l$-periodic solution of (\ref{xis1}) for $\tilde \xi^0_s$.

A choice of a $\l$-periodic $\partial_U$-invariant function $c_s(z,t)$ does not affect the
periodicity property of $\tilde \xi_s=c_s+\tilde\xi_s^0$. It changes the right hand side
of (\ref{ap1}). A particular solution of the new equation is given by the
formula $\chi_s^0=\chi_s^*+c_s\tilde \xi^0_1$. Therefore,
$\tilde \xi_{s+1}^0=\p_t \tilde \xi_s+\chi_s^0$ is a $\l$-periodic solution of (\ref{xis1})
for $\tilde \xi_s$. The choice of $c_s$ does affect the existence
of periodic solutions of
the equation
\beq\label{ap2}
\p_U\chi_{s+1}=-(u+b)\tilde\xi_{s+1}^0-b\p_t\tilde\xi_{s}
\eeq
Let $\tilde\chi_{s+1}$ be a solution of
the equation
\beq\label{ap3}
\p_U\tilde\chi_{s+1}=-(u+b)\tilde\xi_{s+1}^*-b\p_t\tilde\xi^0_{s}\,.
\eeq
Then the function
\beq\label{v7}
\chi_{s+1}(z,t)=\tilde \chi_{s+1}(z,t)+{1\over 2}c_s(z,t)\,(\xi_1^0(z,t))^2
+(\xi_1^0(z,t)-(l,z)b)\p_t c_s (z,t),
\eeq
is a solution of (\ref{ap2}). In order to make  $\chi_{s+1}$ periodic,
the function $c_s(z,t)$ should satisfy the linear differential equation
\beq\label{kp41}
\p_t c_s(z,t)=((l,\l)b)^{-1}(\tilde\chi_{s+1}(z+\l,t)-\tilde\chi_{s+1}(z,t))\,.
\eeq
This equation, together with the initial condition $c_s(z)=c_s(z,0)$ uniquely defines
$c_s(z,t)$.
The induction step is then completed. We have shown that the ratio of two
periodic formal series $\phi_1$ and $\phi$ is $t$-independent. Therefore,
equation (\ref{v2}), where $\rho(z,k)$ is defined by the evaluation of
both the sides at $t=0$, holds. The lemma is thus proven.

\begin{cor} Let $\l_1,\ldots,\l_d$ be a set of linear independent vectors of
the lattice $\Lambda_U$ and let $z_0$ be a point of $\mathbb C^d$. Then, under
the assumptions of the previous lemma,
 there is a unique wave solution of equation (\ref{sh}) such that
the corresponding formal series $\phi(z,t,k;z_0)$ is quasi-periodic with respect
to $\Lambda_U$, i.e. for
$\l\in \Lambda_U$
\beq\label{v10}
\phi(z+\l,t,k;z_0)=\phi(z,t,k;z_0)\,\mu_{\l}(k)
\eeq
and satisfies the normalization conditions
\beq\label{v11}
\mu_{\l_i}(k)=1,\ \ \ \phi(z_0,0,k;z_0)=1.
\eeq
\end{cor}
The proof is identical to that in the part (b) of Lemma 12 in \cite{shiota}
(compare with the proof of the corollary in \cite{kr-schot}).

\section{The spectral curve}

In this section we show that $\l$-periodic wave solutions of equation (\ref{sh}), with
$u$ as in (\ref{anz}), are common eigenfunctions of rings of commuting operators
and identify $X$ with the Prym variety of the spectral curve of these rings.

Note that a simple shift $z\to z+Z$, where $Z\notin \Sigma,$ gives
$\l$-periodic wave solutions with meromorphic coefficients along the affine
subspaces $Z+\mathbb C^d$. These $\lambda$-periodic wave solutions are related to each other
by $t$-independent, $\p_U$-invariant factor. Therefore choosing in the neighborhood
of any $Z\notin \Sigma,$ a hyperplane orthogonal to the vector $U$ and
fixing initial data on this hyperplane at $t=0,$ we define the corresponding
series $\phi(z+Z,t,k)$ as a {\it local} meromorphic function of $Z$ and the
{\it global} meromorphic function of $z$.

\begin{lem} Let equation (\ref{t}) hold. Then there is a unique
pseudo-differential operator
\beq\label{LL}
\L(Z,\p_t)=\p_t+\sum_{s=1}^{\infty} w_s(Z)\p_t^{-s}
\eeq
such that for $Z+Vt\notin \Sigma$
\beq\label{kk}
\L(Ux+Vt+Z,\p_t)\,\psi=k\,\psi\,,
\eeq
where $\psi=e^{(kt+bxk^{-1})} \phi(Ux+Z,t,k)$ is a $\l$-periodic solution of
(\ref{sh}).
The coefficients $w_s(Z)$ of $\L$  are meromorphic functions on the abelian variety $X$ with
poles along  the divisor $\Theta$.
\end{lem}
{\it Proof.} The construction of $\L$ is standard for the KP theory. First, we define
$\L$ as a pseudo-differential operator with coefficients $w_s(Z,t)$,
which are functions of $Z$ and $t$.

Let $\psi$ be a
$\l$-periodic wave solution. The substitution of (\ref{psi2}) into (\ref{kk})
gives a system of equations
that recursively define $w_s(Z,t)$ as differential polynomials in $\tilde\xi_s(Z,t)$.
The coefficients of $\psi$ are local meromorphic functions of $Z$, but
the coefficients of $\L$ are well-defined
{\it global meromorphic functions} of on $\mathbb C^g\setminus\Sigma$, because
different $\l$-periodic wave solutions are related to each other by $t$-independent
factor, which does not affect $\L$. The singular locus is
of codimension $\geq 2$. Then Hartogs' holomorphic extension theorem implies that
$w_s(Z,t)$ can be extended to a global meromorpic function on $\mathbb C^g$.

The translational invariance of $u$ implies the translational invariance of
the $\l$-periodic wave solutions. Indeed, for any constant $s$, the series
$\phi(Vs+Z,t-s,k)$ and $\phi(Z,t,k)$ correspond to $\l$-periodic solutions
of the same equation. Therefore, they coincide up to a $t$-independent,
$\p_U$-invariant factor.
This factor does not affect $\L$. Hence, $w_s(Z,t)=w_s(Vt+Z)$.

The $\l$-periodic wave functions corresponding to $Z$ and
$Z+\lambda'$ for any $\lambda'\in \Lambda$
are also related to each other by a $t$-independent, $\partial_U$-invariant factor.
Hence, $w_s$ are periodic with respect to $\Lambda$ and therefore are
meromorphic functions on the abelian variety $X$.
The lemma is proved.

\begin{lem} Let $\L$ be a pseudo-differential operator corresponding to $\l$-periodic
solution and $\L^*$ be its formal adjoint operator. Then the following equation
\beq\label{ad1}
\L^*=-\p_t\L\p_t^{-1}
\eeq
holds.
\end{lem}
Recall, that the operator which is formally adjoint to $(w\p^{i})$ is the operator
$(-\p)^i\cdot w$, where $w$ stands for the operator of multiplication by the function $w$.
Below we will use the notion of the left action of an operator which is identical to
the formal adjoint action, i.e. by definition we assume that for a function $f$
the identity
\beq\label{ad2}
(f\D)=\D^*f
\eeq
holds.

\noindent
{\it Proof.} If $\psi$ is as in Lemma 3.3, then  there exists a unique pseudo-differential
operator $\Phi$ such that
\beq\label{S}
\psi=\Phi e^{kt},\ \ \Phi=1+\sum_{s=1}^{\infty}\f_s(Ux+Z,t)\p_t^{-s}.
\eeq
The coefficients of $\Phi$ are universal differential polynomials on $\tilde\xi_s$.
Therefore, $\f_s(z+Z,t)$ is a global meromorphic function of $z\in C^d$ and
a local meromorphic function of $Z\notin \Sigma$.
Note that $\L=\Phi(\p_t)\, \Phi^{-1}$, and the equation $H\psi=0$ is equivalent to the
operator equation
\beq\label{ff}
\p_t \cdot\Phi_x+u\Phi=0,
\eeq
where $\Phi_x$ is the pseudo-differential operator with the coefficients $\p_x\f$.
Note that (\ref{ff}) implies
\beq\label{mos1}
\p_x\L=[\p_t^{-1}\cdot u, \L].
\eeq
Let us define the dual  wave function $\psi^*$ by the formula
\beq\label{ad3}
\psi^*=\left(e^{-kt}\p_t\cdot\Phi^{-1}\cdot\p_t^{-1}\right)=
\left(\p_t^{-1}\left(\Phi^{-1}\right)^*\p_t\right)e^{-kt}.
\eeq
Equation (\ref{ff}) implies $H\psi^*=0$. The dual wave function $\psi^*$
is $\l$-periodic. Therefore, the same arguments as used above show that
if equation (\ref{sys}) is satisfied, then the dual wave function is of the form
$\psi^*=e^{-(kt+bxk^{-1})}\phi^*(Ux+Z,t,k)$, where
the coefficients $\tilde\xi_s^*(z+Z,t)$ of the formal series
\beq\label{psi2+}
\phi^*(z+Z,t,k)=1+\sum_{s=1}^{\infty}\tilde\xi^*_s(z+Z,t)\, k^{-s}
\eeq
have simple poles at the divisor $\Theta^{\,U}(t)$. They are
$\l$-periodic. Therefore,

\beq\label{ad4}
\phi^*(z+Z,t,k)=\phi(z+Z,t,-k)\rho(z+Z,k),
\eeq
where $\rho$ is a $t$-independent, $\p_U$-invariant factor. Equation (\ref{ad4}) implies
(\ref{ad1}) and the lemma is proved.

\bigskip
\noindent{\bf Commuting differential operators}.
Let as denote {\it strictly} positive differential part of the pseudo-differential operator $\L^m$ by
$\L^m_+$, i.e. if $\L^m=\sum_{i=-m}^{\infty} F_{m}^{(i)}\p_t^{-i}$, then \footnote{
Note that this definition differs from the one used in the KP theory, where plus subscript
denotes nonnegative part of a pseudo-differential operator}
\beq\label{ad5}
\L^m_+=\sum_{i=1}^{m} F_{m}^{(-i)}\p_t^{i},\ \ \ \
\L^m_-=\L^m-\L^m_+=F_m^{(0)}+F_m^{(1)}\p_t^{-1}+O(\p^{-2}).
\eeq
By definition of the residue of a pseudo-differential operator, the  first leading coefficients
of $\L^m_-$ are
\beq\label{res1}
F_m^{(0)}={\rm res}_{\,\p}\  \left(\L^m\p_t^{-1}\right),\ \ \ F_m^{(1)}={\rm res}_{\,\p}\  \L^m.
\eeq

\begin{lem} The operators $\L^m_+$ satisfy the equations
\begin{eqnarray}
H\L^{2m}_+\ \ &=&-F_{2m,\,x}^{(0)}\,\p_t-{1\over 2} F_{2m,\,xt}^{(0)}+B_{2m}H,
\label{sb3}\\
H\L^{2m+1}_+&=&-F_{2m+1,x}^{(1)}+B_{2m+1}H, \label{sb4}
\end{eqnarray}
where $B_m$ is a pseudo-differential operator in the variable $t$.
\end{lem}
{\it Proof.} First, we prove the equation
\beq\label{ad8}
H\L^m_+=-F_{m,\,x}^{(0)}\,\p_t-\left(F_{m,\,xt}^{(0)}+F_{m,\,x}^{(1)}\right)+B_mH.
\eeq
Each operator $\D$ of the form $\D=\sum_{i=N}^{\infty} (a+b\p_x)\p_t^{-i}$
can be uniquely represented in the form $\D=D_1+D_2H$,
where $D_{1,2}$ are pseudo-differential operators in the variable $t$.
Consider such a representation for the operator $H\L^m=D_1+D_2H$. From the definition
of $\L$ it follows that $H\L^m\psi=0$. That implies $D_1=0$ or the equation
\beq\label{ad12}
H\L^m=D_2H.
\eeq
We have the identity
\beq\label{ad9}
[\p_x\p_t+u,\L^m_+]=\L_{+,\,xt}^m+\L^m_{+,x}\p_t+[u,\L^m_+]-\L^m_{+,\,t}\p_t^{-1}\cdot u
+\L^m_{+,\,t}\cdot\p_t^{-1}\cdot H\,.
\eeq
The first three terms are differential operators in the $t$ variable. By definition of
$\L^m_+$ the fourth term is also a differential operator.
Therefore, the pseudo-differential operator $D_1$ in the decomposition
$H\L^m_+=D_{m,1}+B_mH$ is a {\it differential} operator.

In the same way we get the equation
\beq\label{ad13}
H\L^m_-=\tilde D_{m,1}+\tilde B_{2m}H,
\eeq
where
\beq\label{ad14}
\tilde D_{m,1}=\L_{-,\,xt}^m+\L^m_{-,x}\p_t+[u,\L^m_-]-\L^m_{-,\,t}\p_t^{-1}\cdot u
\eeq
By definition of $\L^m_-$ the operator $\tilde D_{m,1}$ is a pseudo-differential
operator of order not greater than $1$.
Equation (\ref{ad12}) implies $H\L_+^m=-H\L^m_-+D_2H$. Hence,
$D_{m,1}=-\tilde D_{m,1}$ is a differential operator of the order
$1$, i.e. has the form $a\p_t+b$. The coefficients of this operator can be easily
found from the leading coefficients of the right hand side of (\ref{ad14}).
Direct computations give equation (\ref{ad8}).

Now in order to complete the proof of (\ref{sb3}) and (\ref{sb4}) it is enough
to use (\ref{ad1}). From equation (\ref{ad1}) and the relation
${\rm res}_{\p} D=-{\rm res}_{\p} D^*$ it follows that
\beq\label{sb1}
F_{2m}^{(1)}=-{\rm res}_{\,\p}\  (\L^*)^{2m}=
-{\rm res}_{\,\p}\ \left(\p_t\L^m\p_t^{-1}\right)=-F_{2m, t}^{(0)}-F_{2m}^{(1)}.
\eeq
In the same way we get
\beq\label{sb5}
F_{2m+1}^{(0)}={\rm res}_{\,\p}\  (\L^{2m+1}\p_t^{-1})=
{\rm res}_{\,\p}\ \left(\L^{2m+1}\p_t^{-1}\right)^*=-F_{2m+1}^{(0)}=0.
\eeq
Equations (\ref{ad8}, \ref{sb1}, \ref{sb5}) imply (\ref{sb3}) and (\ref{sb4}). The lemma
is proved.

The following statement is a direct corollary of equations (\ref{sb1}, \ref{sb5}).
\begin{cor} The operators $\L^m_+$ satisfy the relation
\beq\label{p}
\left(\L^m_+\right)^*=(-1)^m \p_t \cdot\L^m_+\cdot\p_t^{-1}
\eeq
\end{cor}
The next step is crucial for the construction of commuting operators.
\begin{lem} The functions $F_{2m}^{(0)}, F_{2m+1}^{(1)}$ have at most second order pole
on the divisor $\Theta$.
\end{lem}
{\it Proof.} The ambiguity in the definition of $\psi$ does not affect the product
\beq\label{J0}
\psi^*\psi=\left(e^{-kt}\p_t\Phi^{-1}\,\p_t^{-1}\right)\left(\Phi e^{kt}\right).
\eeq
Therefore, although each factor is only a local meromorphic function on
$\mathbb C^g\setminus \Sigma$, the coefficients $J_s^{(0)}$ of the product
\beq\label{J}
\psi^*\psi=\phi^*(Z,t,k)\,\phi(Z,t,k)=1+\sum_{s=2}^{\infty}J_s^{(0)}(Z,t)k^{-s}
\eeq
are {\it global meromorphic functions} of $Z$. Moreover, the translational invariance
of $u$ implies that they have the form $J_s(Z,t)=J_s(Z+Vt)$.
Each of the factors in the left hand
side of (\ref{J}) has a simple pole on $\Theta-Vt$.
Hence, $J_s(Z)$ is a meromorphic function
on $X$ with a second order pole at $\Theta$.

From the definition of $\L$ it follows that
\beq\label{20}
\res_k\left(\psi^*(\L^n\psi)\right){k^{-1}dk}=
\res_k\left(\psi^*k^n\psi\right){k^{-1}dk}=J_{n}^{(0)}.
\eeq
On the other hand, using the identity
\beq\label{dic}
\res_k \left(e^{-kx}\D_1\right)\left(\D_2e^{kx}\right)dk=\res_{\p}\left(\D_2\D_1\right),
\eeq
we get
\beq\label{201}
\res_k\left(\psi^*\left(\L^n\psi\right)\right)\,k^{-1}dk=
\res_k\left(e^{-kt}\p_t\Phi^{-1}\p_t^{-1}\right)\left(\L^n\Phi \,\p_t^{-1}e^{kt}\right)dk=
\res_{\p}\left(\L^n\p_t^{-1}\right)=F_n^{(0)}.
\eeq
Therefore, $F_n^{(0)}=J_{n}^{(0)}$ has the second order pole at $\Theta$.

Consider now the coefficients $J_s^{(1)}$ of the series
\beq\label{J1}
\psi^*\psi_t-\psi_t^*\psi=2k+\sum_{s=1}^{\infty}2J_s^{(1)}(Z,t)k^{-s}.
\eeq
They are meromorphic functions on $X$ with the second order pole at $\Theta$. We have
\beq\label{sb6}
2J_n^{(1)}={\rm res}_k \left((\psi^*\L^n)\psi_t-\psi_t^*(\L^n\psi)\right)k^{-1}dk=
{\rm res}_{\p} (\L^n+\p_t\L^n\p_t^{-1})=2F_n^{(1)}+F_{n,t}^{(0)}
\eeq
Then from equation (\ref{sb5}) it follows that $F_{2m+1}^{(1)}=J_{2m+1}^{(1)}$
and the lemma is proved.

Let ${\bf F}$ be a direct sum of the linear spaces
${\bf {\hat F}^{\a}}, \ \a=0,1,$ spanned by
$\{F_{2m+\a}^{(\a)}, \ m=0,1,\ldots\}$. They are subspaces of
the $2^g$-dimensional space of the abelian functions with at most second order
pole at $\Theta$. Therefore, for all but $\hat g^{\a}=\dim\ {\bf \hat F}^{\a}$ positive
integers $2n+\a$, there exist constants $c_{i,n}^{(\a)}$ such that
\beq\label{f1}
F_{2n+\a}^{(\a)}(Z)+\sum_{i=1}^{n} c_{i,n}^{(\a)}F_{2n-2i+\a}^{(\a)}(Z)=0.
\eeq
Let $I^{(\a)}$ denote the subset of integers $2n+\a$
for which none of such constants exist.
We call the union $I=I^{(0)}\cup I^{(1)}$ the gap sequence.
\begin{lem} Let $\L$ be the pseudo-differential operator corresponding to
a $\l$-periodic wave function $\psi$ constructed above.
Then, for the differential operators
\beq\label{a2}
L_{2n+\a}=\L^{2n+\a}_++\sum_{i=1}^{n} c_{i,n}^{(\a)}\L^{2n+\a-2i}_+, \ \ \ 2n+\a\notin I^{\a},
\eeq
the equations
\beq\label{lp}
L_{2n+\a}\,\psi=a_{2n+\a}(k)\,\psi, \ \ \
a_{2n+\a}(k)=k^{2n+\a}+\sum_{s=1}^{\infty}a_{s,n}k^{2n+\a-s},
\eeq
where $a_{s,n}$ are constants, hold.
\end{lem}
{\it Proof.} First note, that from (\ref{sb3}, \ref{sb4}) it follows that
$HL_{2n+\a}\psi=0$.
Hence, if $\psi$ is a $\l$-periodic wave solution of (\ref{sh})
corresponding to $Z\notin \Sigma$, then $L_{2n+\a}\psi$ is also a formal
solution of the same equation. That implies the equation
$L_{2n+\a}\psi=a_{2n+\a}(Z,k)\psi$, where $a_{2n+\a}(Z,k)$ is $t$-independent and
$\p_U$-invariant function of the variable $Z$.
The ambiguity in the definition of
$\psi$ does not affect $a_{2n+\a}$. Therefore, the coefficients of $a_{2n+\a}$
are well-defined
{\it global} meromorphic functions on $\mathbb C^g\setminus \Sigma$. The $\p_U$-
invariance of $a_{2n+\a}$ implies that $a_{2n+\a}$, as a function of $Z$, is holomorphic
outside the locus. Hence it has an extension to a holomorphic function on $\mathbb C^g$.
The $\l$-periodic wave functions corresponding to $Z$ and
$Z+\lambda'$ for any $\lambda'\in \Lambda$ are related to each other by a $t$-independent,
$\partial_U$-invariant factor.
Hence, $a_{2n+\a} $ is periodic with respect to $\Lambda$ and therefore is $Z$-independent.
Note that $a_{2s+1,n}=0$ and $a_{2s,n}=c_{s,n}$ if $s\leq n$.
The lemma is proved.

The operator $L_m$ can be regarded as a $(Z,x)$-parametric family   of
ordinary differential operators $L_m^Z$ whose coefficients have the form
\beq\label{lu}
L_m^{Z,x}=\p_t^m+\sum_{i=1}^m u_{i,m}(Ux+Vt+Z)\, \p_t^{m-i},\ \ m\notin I.
\eeq
where $u_{i,m}(Z)$ are abelian function regular outside of $\Theta$.
For $Z+Ux\notin \Sigma_-$ the coefficients of $L_m^{Z,x}$ are meromorphic functions
of the variable $t$, which are not identically equal infinity. Recall, that
$\Sigma_-$ is a $\p_V$-invariant set of $\Theta$.
\begin{cor} The operators $L_m^{Z,x}$
commute with each other,
\beq\label{com1}
[L_n^{Z,x},L_m^{Z,x}]=0.
\eeq
\end{cor}
From (\ref{lp}) it follows that $[L_n^{Z,x},L_m^{Z,x}]\psi=0$. The commutator is
an ordinary differential operator. Hence, the last equation implies (\ref{com1}).

\begin{lem} Let $\A^{Z,\,x},\ Z+Ux\notin \Sigma_-,$ be a commutative ring of ordinary
differential operators spanned by the operators $L_n^{Z,x}$.
Then there is an irreducible algebraic
curve $\G$ such that $\A^{Z,\,x}$ is isomorphic
to the ring $A_-(\G,P_+,P_-)$ of the meromorphic functions on $\G$ with
the only pole at a smooth point $P_-$ vanishing at another smooth point
$P_+$. The correspondence $Z\to \A^{Z,\,0}$ defines a holomorphic
map of $X\setminus \Sigma_-$ into the space of torsion-free rank 1 sheaves
$\F$ on $\G$
\beq\label{is}
j: X\backslash\Sigma_-\longmapsto \overline{\rm Pic}(\G).
\eeq
On an open set the map $j$ is an imbedding.
\end{lem}
The proof of the lemma is almost identical to
the proof of lemma 3.4 in \cite{kr-schot}.
It is the fundamental fact of the theory of commuting linear
ordinary differential operators (\cite{kr1,kr2,ch1,ch2,mum}) that there is a
natural correspondence
\beq\label{corr}
\A\longleftrightarrow \{\G,P_-, [k^{-1}]_1, \F\}
\eeq
between {\it regular} at $t=0$ commutative rings $\A$ of ordinary linear
differential operators in the variable $t$,
containing a pair of monic operators of co-prime orders, and
sets of algebraic-geometrical data $\{\G,P_-, [k^{-1}]_1, \F\}$, where $\G$ is an
algebraic curve with a fixed
first jet $[k^{-1}]_1$ of a local coordinate $k^{-1}$ in the neighborhood of a smooth
point $P_-\in\G$ and $\F$ is a torsion-free rank 1 sheaf on $\G$ such that
\beq\label{sheaf}
H^0(\G,\F)=H^1(\G,\F)=0.
\eeq
The correspondence becomes one-to-one if the rings $\A$ are considered modulo conjugation
$\A'=g(t)\A g^{-1}(t)$.

Note, that in \cite{kr1,kr2,ch1,ch2} the main attention was paid to the generic case of
the commutative rings corresponding to smooth algebraic curves.
The invariant formulation of the correspondence given above is due to Mumford \cite{mum}.

The algebraic curve $\G$ is called the spectral curve of $\A$.
The ring $\A$ is isomorphic to the ring $A(\G,P_-)$ of meromorphic functions
on $\G$ with the only pole at the puncture $P_-$. The isomorphism is defined by
the equation
\beq\label{z2}
L_a\psi_0=a\psi_0, \ \ L_a\in \A, \ a\in A(\G,P_-).
\eeq
Here $\psi_0$ is a common eigenfunction of the commuting operators. At $t=0$ it is
a section of the sheaf $\F\otimes\O(-P_-)$.

{\bf Important remark}. The construction of the correspondence (\ref{corr})
depends on a choice of the initial point $t_0=0$. The spectral curve and the sheaf $\F$
are defined by the evaluations of the coefficients of generators of $\A$ and a finite
number of their derivatives at the initial point. In fact, the spectral curve
is independent on the choice of $t_0$, but the sheaf does depend on it, i.e. $\F=\F_{t_0}$.

Using the shift of the initial point it is easy to show that the correspondence
(\ref{corr}) extends to the commutative rings of operators whose coefficients are
{\it meromorphic} functions of $t$ at $t=0$. The rings of operators having poles at $t=0$
correspond to sheaves for which the condition (\ref{sheaf}) is violated.

As it was mentioned above, the operators $L_n$, $L_m$ can be seen as a $(Z,x)$-parametric family
of commuting ordinary differential operators in the variable $t$. Let $\G^{Z,x}$ be
the corresponding spectral curve. The eigenvalues $a_n(k)$ of the operators $L_n^{Z,x}$
defined in (\ref{lp}) coincide with the Laurent expansions at $P_-$ of the meromorphic
functions $a_n\in A(\G^{Z,x},P_-)$. They are $(Z,x)$-independent. Hence, the spectral curve
is $(Z,x)$-independent, as well, $\G=\G^{Z,x}$.

Equations (\ref{t}), which are equivalent
to (\ref{sys}) and are sufficient for the construction of the $\l$-periodic wave solutions,
are symmetric with respect to $x$ and $t$. Therefore, the simple interchange of the
variables $x$ and $t$ shows that if equations (\ref{t}) hold then there exist
commuting ordinary differential operators $L_m^+$
of the form
\beq\label{lu1}
L_m^+=\p_x^m+\sum_{i=1}^m u_{i,m}^+(Ux+Vt+Z)\, \p_x^{m-i},\ \ m\notin I^+,
\eeq
where $I^+$ is the gap sequence associated with the variable $x$.
These operators satisfy the equations
\beq\label{10+}
L^+_mH=B_m^+H,
\eeq
where $B_m^+$ are differential operators in the variable $x$.

Let $\psi$ be a $\l$-periodic solution of (\ref{sh}).
Then the same arguments as in the proof of Lemma 4.5 show that equations
(\ref{10+}) imply
\beq\label{lp+}
L_{n}^+\,\psi=a_{n}^+(k)\,\psi, \ \ \
a_n^+=\sum_{s=1}^{\infty}a^+_{s,n}k^{-s}\,,
\eeq
where $a_{s,n}^+$ are constants. From (\ref{lp+}) it follows that the operators
$L_n,L_m^+$ satisfy the equation
\beq\label{200}
[L_n,L_m^+]=BH,
\eeq
where $B$ is a differential operator in the variables $x,t$.
Equation (\ref{lp+}) also implies that there exists a polynomial $\tilde R$ such that
$\tilde R(L_n,L_m^+)\psi=0$, i.e. eigenvalues $a_n, a_m^+$ of $L_n$ and $L_m^+$
satisfy the equation $\tilde R(a_n,a_m^+)=0$. Therefore, the spectral curves of commutative
rings $\A^Z$ and $\A^Z_+$ coincide. Note, that ((\ref{lp+}) implies that $a_m^+$ vanishes
at $P_-$. The symmetry between $x$ and $t$ variables implies that
the ring $\A_+^Z$ is isomorphic to the ring $A_+(\G,P_+,P_-)$ of meromorphic functions on
$\G$ with the only pole at $P_+$ that vanish at $P_-$.

Let us fix $x=0$ and consider the commuting operators $L_n^Z=L_n^{Z,0}$.
The construction of the correspondence
(\ref{corr}) implies that if the coefficients of the operators in $\A$
holomorphically depend on parameters, then the algebraic-geometrical spectral data are
also holomorphic functions of the parameters.

Therefore,  $j$ is holomorphic out of $\Theta$.
Then, using the shift of the initial point and the fact, that $\F_{t_0}$ holomorphically
depends on $t_0$, we get that $j$ holomorphically extends on $\Theta\setminus \Sigma_-$,
as well.

The theta-divisor is not invariant under the shifts by constant vectors.
Hence, for the generic $Z$ and $Z'$ the operators in $\A^Z$ and $\A^{Z'}$ have
different poles. Hence, those rings do not coincide. Thus, the map $j$ is
an imbedding on an open set. The lemma is proved.

It implies the global existence of the wave function.
\begin{lem} Let equations (\ref{t}) hold. Then there exists a common
eigenfunction of the operators $L_n^Z$ of the form
$\psi=e^{kt}\phi(Vt+Z,k)$ such that
the coefficients of the formal series
\beq\label{psi6}
\phi(Z,k)=1+\sum_{s=1}^{\infty}\xi_s(Z)\, k^{-s}
\eeq
are global meromorphic functions with a simple pole at $\Theta$.
\end{lem}
The proof of the lemma is identical to the proof of lemma 3.5 in \cite{kr-schot}.
The function $\psi$ is first defined for $Z\notin \Sigma_-$ as the inverse
image $\psi=j^*\hat\psi_{BA}$ of the Baker-Akhiezer function, which is known to be
globally defined on $\overline {\rm Pic}(\G)$. Then, Hartogs' extension theorem
implies that $\psi$ has a meromorphic extension on $C^g$. The Baker-Akhiezer
function is regular out of divisor
corresponding to the commutative rings of operators whose coefficients have
poles at $t=0$. Hence, $\psi$ is holomorphic out of $\Theta$.

Let us show now that the correspondence $\psi\to \psi^*$ defines
an involution of $\G$ under which $P_{\pm}$ are fixed.

\begin{lem} The eigenvalues $a_{2n+\a}$ of the commuting operators $L_{2n+\a}$ satisfy
the relation
\beq\label{aa3}
a_{2n+\a}(k)=(-1)^\a a_{2n+\a}(-k).
\eeq
\end{lem}
{\it Proof.} From equations (\ref{ad4},\ref{p},\ref{a2}) it follows that
\beq\label{aa4}
\psi^*_t(L_{2n+\a}\psi)=a_{2n+\a}(k)(\psi^*_t\psi),
\eeq
\beq\label{aa5}
(\psi^*_tL_{2n+\a})\psi=\left((L_{2n+\a}^*\psi^*_t)\psi\right)=(-1)^{\a}
a_{2n+\a}(-k)(\psi^*_t\psi).
\eeq
The left and  right action of pseudo-differential operators are formally adjoint,
i.e., for any two operators the equality $\left(e^{-kt}\D_1\right)\left(\D_2e^{kt}\right)=
e^{-kt}\left(\D_1\D_2e^{kt}\right)+\p_t\left(e^{-kt}\left(\D_3e^{kt}\right)\right)$
holds. Here $\D_3$ is a pseudo-differential operator whose coefficients are differential
polynomials in the coefficients of $\D_1$ and $\D_2$.
Therefore, equations (\ref{aa4},\ref{aa5}) imply
\beq\label{aa6}
\left(a_{2n+\a}(k)-(-1)^{\a}a_{2n+\a}(-k)\right)(\psi^*_t\psi)=\p_t Q_{2n+\a}.
\eeq
The coefficients of the series $Q_{2n+\a}$ are differential polynomials
on the coefficients of the wave operator $\Phi$ defined by equation (\ref{S}).
For the globally defined wave function $\psi$,
which exists according to the previous lemma, the coefficients of the wave operator are
global meromorphic functions.
Hence,
\beq\label{ap5}
Q_{2n+\a}=\sum_{s=1}^{\infty}Q_{2n+\a,s}(Vt+Z),
\eeq
where $Q_{2n+\a,s}(Z)$ are  meromorphic functions regular out of $\Theta$.

In a similar way we have
\beq\label{aa7}
\psi^*_t\psi=\left(e^{-kt}\p_t\Phi^{-1}\right)\left(\Phi e^{kt}\right)=
k+\p_t Q^{(1)}.
\eeq
The series $Q^{(1)}$ has the form
\beq\label{ap6}
Q^{(1)}=\sum_{s=1}^{\infty}Q_s^{(1)}(Vt+Z),
\eeq
where $Q_s^{(1)}(Z)$ are  meromorphic  functions regular out of $\Theta$.

Let us fix a neighborhood of the theta-divisor in $X$.
It defines the neighborhood $S$ of $\Theta$ in ${\mathbb C}^{\wh g}$.
Outside of $S$ the functions $Q_s^{(1)}$ are bounded.
Consider a sequence of real numbers $l_i\to \infty$ such
that $Z\pm Vl_i$ is not in $S$. Then, from (\ref{aa7}) it follows that
\beq\label{aa9}
\langle \psi^*_t\psi\rangle=
\lim_{l_i\to \infty} {1\over 2l_i}\int_{-l_i}^{l_i} (\psi^*_t\psi)dt=k\,.
\eeq
The integration (\ref{aa9}) is taken along a curve connecting points
$Z+\pm Vl_i$ and which does not intersect $\Theta$.

The same arguments imply that under "averaging" in $t$ the right hand side of (\ref{aa6})
vanishes. Hence, (\ref{aa6}) implies (\ref{aa3}). The lemma is proved.

The series $a_n(k)$ are the expansions at $P_-$ of meromorphic functions
on $\G$. Therefore, from (\ref{aa3}) it follows there exists a holomorphic involution
$\s:\G\to\G$ of the spectral curve such that
\beq\label{ap7}
a_n^{\s}=a_n(\s(P))=(-1)^na_n(P)
\eeq
The point $P_-$ is fixed under $\s$ and the local parameter is odd with respect
to $\s$, i.e. $\s^*k=-k$. In the same way using $x$ variable instead of $t$ we get
that the second puncture $P_+\in \G$ is also fixed under $\s$.

The involution $\s$ induces an involution on the generalized Jacobian
$J(\G)$ which is by definition is the group of the equivalence classes
of zero-degree divisors on $\G$, i.e. $J(\G)={\rm Pic^0}(\G)$.
The odd subgroup of $J(\G)$ with respect to the induced involution $\s^*$ is
the Prym variety of the spectral curve, $\P(\G)={\rm ker} (1+\s^*)$.
Our next goal is to show that $\P(\G)$ of the spectral curve
is {\it compact}.
\begin{lem} There exist $g$-dimensional vectors $V_{2m+1}=\{V_{2m+1,k}\}$
and constants $v_{2m+1}$ such that
\beq\label{nnov7}
F_{2m+1}^{(1)}(Z)=\sum_{k=1}^g V_{2m+1,k} \p_Vh_k(Z)+v_{2m+1},
\eeq
where $F_{2m+1}^{(1)}=\res_{\p}\, \L^{2m+1}$ and
$h_k=\p_{z_k}\ln \theta(Z)$.
\end{lem}
{\it Proof.} From equations (\ref{J1}, \ref{sb6}) and (\ref{aa7}) it follows that
$F_{2m+1}^{(1)}=-2\p_V Q_{2m+1}^{(1)}$, where $Q_{2m+1}^{(1)}$ is a
meromorphic function with a pole along $\Theta$. The function
$F_{2m+1}^{(1)}$ is an abelian function. Hence, for any vector $\l$ in the period
lattice $Q_{2m+1}^{(1)}(Z+\l)=Q_{2m+1}^{(1)}(Z)+c_{m,\l}$.
There is no abelian function with a simple pole on $\Theta$. Hence, there exists
a constant $q_n$ and $g$-dimensional vectors $l_m$ and $V_{2m+1}$, such that
$Q_{2m+1}^{(1)}=q_m+(l_m,Z)+(V_{2m+1},h(Z))$, where $h(Z)$ is a vector with the
coordinates
$h_k(Z)$. Therefore, $F_{2m+1}^{(1)}=(l_m,V)+(V_{2m+1},\p_V h(Z))$.
The lemma is proved.

In order to complete the proof of our main result we need few more facts
of the KP theory: flows of the KP hierarchy define deformations of the commutative
rings $\A$ of ordinary linear differential operators. The spectral curve
is invariant under these flows. If a commutative ring $\A$ of linear ordinary
differential operators is {\it maximal}, i.e., it is not contained in any bigger
commutative ring, then the KP orbit of $\A$ is isomorphic to the generalized Jacobian
$J(\G)={\rm Pic}^0 (\G)$ of the spectral curve of $\A$
(see details in \cite{shiota,kr1,kr2,wilson}).

The KP hierarchy in the Sato form is a system of  commuting differential
equation for a pseudo-differential operator $\L$
\beq\label{z4}
\p_{t_n}\L=[\L^n_+,\L]\,.
\eeq
If the operator $\L$ is as above. i.e., if it is defined by $\l$-periodic wave
solutions of equation (\ref{sh}), then equation (\ref{sb4}) implies
that for odd $n$ equations (\ref{z4}) are equivalent to
the equations
\beq\label{z5}
\p_{t_{2n+1}}u=\p_x F_{2m+1}^{(1)}(Ux+Vt+Z).
\eeq
The first time of the hierarchy is identified with the variable $t_1=t$.

Equations (\ref{z5}) identify the space generated by the functions
$\p_U F_{2m+1}^{(1)}$ with the tangent space at  $\A^Z$ of the orbit of the part
of the NV hierarchy associated with the puncture $P_-$.
In terms of $u$ the deformation with respect to $z_i$ is given by the equation
\beq\label{z61}
\p_{z_i}u=\p_x \p_V h_i,\ \ h_i=\p_{z_i}\ln \theta(Z).
\eeq
Equations (\ref{z5}, \ref{z61}) and (\ref{nnov7}) imply
\beq\label{nov}
\p_{t_{2n+1}}=\p_{V_{2n+1}}=\sum_{k=1}^g V_{2n+1,k}\p_{z_k}.
\eeq
Hence, the orbit of $\A^Z$ is isomorphic to the factor of $Z+Y/T(Z)$
of the affine subvariety $Z+Y\subset X$, where $Y$ is the closure in $X$ of the subgroup
$\langle\sum_n V_{2n+1}t_{2n+1}\rangle$, and $T(Z)$ is a lattice in the
universal cover of $Y$.

\begin{lem} For the generic $Z\notin \Sigma_-$, the orbit of $\A^Z$ under
the NV flows defines an isomorphism:
\beq\label{imb}
i_Z:\P(\G)\longmapsto Z+Y\subset X.
\eeq
\end{lem}
{\it Proof.} Recall, that according to \cite{shiota2} the NV orbit of a maximal
commutative ring is isomorphic to the Prym variety of the corresponding spectral curve.
The arguments showing  that $\A^Z$ is maximal for the generic $Z$ are identical to
those used in \cite{kr-schot}.
Indeed, suppose that $\A^Z$ is not maximal for all $Z$.
Then there exits  $2n+\a\in I$, where $I$ is the gap sequence defined above, such
that for each $Z\notin \Sigma_-$ there exists an operator $L_{2n+\a}^Z$
of order $2n+\a$ which commutes with all the operators $L_m^Z\in \A^Z$.
Therefore, it commutes with $\L$.
That implies the equality
\beq\label{z3}
F_{2n+\a}^{(\a)}(Z)+\sum_{i=1}^{n} c_{i,n}^{(\a)}(Z)F_{2n-2i+\a}^{(\a)}(Z)=0.
\eeq
Note the difference between (\ref{f1}) and (\ref{z3}). In the first equation the
coefficients $c_{i,n}^{(\a)}$ are constants.

The $\l$-periodic wave solution of equation
(\ref{sh}) is a common eigenfunction of all  commuting operators, i.e.
$L_{2n+\a}\psi=a_{2n+\a}(Z,k)\psi$, where is $\p_V$-invariant.
The compactness of $X$ implies that $a_{2n+\a}$ is $Z$-independent. The first $n$
coefficients of $a_{2n+\a}$ coincide with the coefficients in (\ref{z3}).
Hence, these coefficients are $Z$-independent. That contradicts the assumption that
$2n+\a \in I$.

The map $j$ defined in Lemma 4.6 restricted to $Z+Y\subset X$ is inverse to $i_Z$.
For the generic $Z$ it is an imbedding. Hence for the generic $Z$
the lattice $T(Z)$ is trivial. The lemma is thus proven.

\begin{cor}
The Prym variety $\P(\G)$ of the spectral curve $\G$ is compact.
\end{cor}
The compactness of the Prym variety is not as restrictive,
as the compactness of the Jacobian (see \cite{mdl}).
Nevertheless, it implies an explicit description of the
the singular points of the spectral curve.
The proof of the following statement is due to Robert Friedman and
is presented in the Appendix.
\begin{cor} (R. Friedman)
The spectral curve $\G$ is smooth outside of fixed points $P_{\pm}, Q_k$ of
the involution $\s$. The branches of $\G$ at $Q_k$ are linear and are not permuted by $\s$.
\end{cor}
An equivalent formulation of the corollary is as follows:
there is a a smooth algebraic curve $\tilde \G$ with involution
$\tilde \s$ and a regular equivariant  map $p:\tilde \G\to \G$ which is one-to-one out
of preimages $Q_k^i,\ i=1\ldots, \nu_k,$ on $\tilde \G$ of the singular points
$Q_k$.

The common eigenfunction of commuting differential operators is well-defined
up to a constant factor for all smooth points of the spectral curve.
It can be analytically extended along the branches of the spectral curve
passing through the singular points, i.e. the preimage $\tilde \psi$
of the Baker-Akhiezer on $\tilde \G$ can be regarded as a section of
a line bundle on $\tilde \G$. From the construction of the correspondence
(\ref{corr}) it follows that the evaluations of $\tilde \psi$ at the preimages
of the singular points $Q_k$ satisfy linear relations
\beq\label{dec1}
\sum_{i=1}^{\nu_k} c_{k,j}^{\,i}\tilde\psi(t,Q_k^i)=0, \ \ j=1,\ldots,n_k.
\eeq
The coefficients of these relations and the zero divisor $D$ of
$\tilde\psi(0,\tilde P)$ can be regarded the data defining the corresponding sheaf $\F$.
The divisor $D$ is the pole divisor of the normalized eigenfunction
$\tilde \psi_0(t,\tilde P)=\tilde \psi(t,\tilde P)/\tilde \psi(0,\tilde P)$.

The following theta-functional formula (\ref{nnov5})
for $\tilde \psi$ is crucial for the final steps
of the proof. First note that
using the transformation $\psi\longmapsto e^{(l(k),Z)}\psi$, where $l(k)$ is a series such
that $(l(k),V)=0$, we may assume without loss of generality that
the series $\phi$ in (\ref{psi6}) satisfies the following monodromy
properties:
\beq\label{nnov2}
\phi(Z+e_j,k)=\phi(Z),\ \ \phi(Z+B_j)=\phi(Z)\rho_j(k)
\eeq
where $e_j$ are the basis vectors in $\mathbb C^g$ and $B_j$ are vectors
defined by the columns of the matrix $B$, corresponding to the principle
polarization of $X$.

Equations (\ref{nnov2}) and the fact the the coefficients of $\phi$ are
meromorphic functions with simple poles along $\Theta$ imply that
there is a series $A(k)$ such that
\beq\label{nnov3}
\phi={\theta(A(k)+Z)\over\theta(Z)}
\eeq
The series $A(k)$ defines an imbedding of the neighborhood of $P_-$ into $X$.

The same arguments  show that there is  a holomorphic map
\beq\label{nnov4}
\tilde A:\tilde \G\longmapsto X
\eeq
such that the function $\tilde\psi(t,\tilde P), \, \tilde P\in \tilde \G,$
given by the formula
\beq\label{nnov5}
\tilde\psi={\theta(\tilde A(\tilde P)+Vt+Z)\over\theta(Vt+Z)}\,e^{t\,\Omega(\tilde P)},
\eeq
is the common eigenfunction of the operators in $\A^Z$. Here
$\Omega(\tilde P)$ is an abelian integral on $\tilde \G$ having the form
$\Omega=k+O(k^{-1})$ at $P_-$. Then the normalized
eigenfunction of the commuting operators is given by the formula
\beq\label{nnov51}
\tilde\psi_0={\theta(\tilde A(\tilde P)+Vt+Z)\,\theta(Z)\over
\theta (\tilde A(\tilde P)+Z) \,\theta(Vt+Z)}\,e^{t\,\Omega(\tilde P)}.
\eeq
Our next goal is to show that the pole divisor
$D$ of $\tilde \psi_0$ satisfies the condition analogous to (\ref{p1}) found
by Novikov and Veselov in the case of the smooth spectral curves.

\begin{lem} The equivalence class of $[D]\in J(\tilde \G)$
of the divisor $D$ satisfies the equation
\beq\label{nnovp}
[D]+[\tilde\s(D)]=K+P_++P_-+\sum_{k,i}Q_k^{\,i}\in J(\tilde\G),
\eeq
where $K$ is the canonical class, i.e. the equivalence class of the
zero-divisor of a holomorphic differential on $\tilde\G$.
\end{lem}
{\it Proof.} Equation (\ref{nnovp}) is equivalent to the condition that the divisor
$D+\s(D)$ is the zero divisor of a meromorphic differential $d\Omega$ on $\tilde\G$
with simple poles at the punctures $P_{\pm}$ and the points $Q_k^i$.
The differential $d\Omega$ is even with respect to the involution and descends to
a meromorphic differential on the
factor-curve $\G_0$.

The existence of such differential can be proved almost identically to the proof
of the statement that the conditions (\ref{p1}) are necessary conditions for
the potential reduction of the $2D$ Schr\"odinger operators given in
\cite{kr-spec} (Theorem 3.1).

Let $\tilde \psi_0(x,t,P)$ be the normalized solution of the Schr\"odinger operator.
It is obtained by the deformation along the $x$-flow from the normalized eigenfunction
of the operators in $\A^Z$ considered above. Therefore, it has the form
(\ref{nv1}) with $\theta_{pr}$ and $A^{pr}$ replaced by $\theta$ and $\tilde A$,
respectively. Following \cite{kr-whith} we present another {\it real} form of
$\tilde\psi_0$. Let us introduce real coordinates of a complex vector
$Z\in {\mathbb C}^{g}$ by the formula $Z=\z'+B\z''$, where $\z',\z''$ are
$g$-dimensional real
vectors, and $B$ is the matrix of $b$-periods of the normalized
holomorphic differentials on $\G$. Then the absolute value $|\phi|$ of the function
$\phi(\z,P), \z=(\z',\z'')$ given by the formula
\beq\label{aa}
\phi(\z,P)={\theta(\tilde A(P)+Z)\over \theta(Z)}e^{2\pi i (\tilde A(P), \,\z'')}
\eeq
is a {\it periodic} function of the coordinates $\z'_k,\z''_k$. For real $x,t$ the function
$\psi$ can be represented in the form
\beq\label{aa8}
\tilde \psi_0={\phi (\wh Ux+\wh Vt+\z,P)\over\phi(\z,P)} e^{tp_-+xp+}.
\eeq
Here $\wh U=(U_1^{+'},U_1^{+''}), \wh V=(U_1^{-'},U_1^{-''})$ are
$2g$-dimensional real vectors
corresponding to the complex vectors $U, V$;
$dp_{\pm}$ are meromorphic differentials on $\tilde \G$ with poles of the second order
at $P_{\pm}$ and whose periods are {\it pure imaginary}.

The differential $d\tilde\psi_0$ is also a solution of the same Schr\"odinger
equation. That implies the equality
\beq\label{nnov31}
\p_x(\p_t\tilde\psi_0^{\s}d\tilde\psi_0-\tilde\psi_0^{\s}\p_td\tilde\psi_0)=
\p_t(\p_x\tilde\psi_0^{\s}d\tilde\psi_0-\tilde\psi_0^{\s}\p_xd\tilde\psi_0)
\eeq
The "averaging" of this equation in the variables $x,t$ gives the equation
\beq\label{nnov32}
\langle\p_t\tilde\psi_{0}^{\s}\tilde\psi_0-\tilde\psi_0^{\s}\p_t\tilde\psi_0\rangle_xdp_-=
\langle\p_x\tilde\psi_0^{\s}\tilde\psi_0-\tilde\psi_0^{\s}\p_x\tilde\psi_0\rangle_t\, dp_+.
\eeq
Here $\langle\cdot\rangle_t$ stands for the mean value in $t$ defined as in (\ref{aa9}),
and $\langle\cdot\rangle_x$ stands for the mean value in $x$ defined in a similar
way. The same arguments as in \cite{kr-spec}, show that the differential
\beq\label{nnov33}
d\Omega={dp_+\over
\langle\p_t\tilde\psi_{0}^{\s}\tilde\psi_0-\tilde\psi_0^{\s}\p_t\tilde\psi_0\rangle_x}=
{dp_-\over
\langle\p_x\tilde\psi_0^{\s}\tilde\psi_0-\tilde\psi_0^{\s}\p_x\tilde\psi_0\rangle_t}
\eeq
is holomorphic on $\tilde \G$ except at the branch points where it has
simple poles. It has zeros at the poles of $\tilde\psi_0$ and $\tilde\psi_0^{\s}$.
The lemma is proved.

The differential $\tilde \psi_0\tilde \psi_0^{\s}d\Omega$ is a meromorphic differential
on $\tilde\G$. Its residues at the points $P_{\pm}$ are equal to $\pm 1$, respectively.
Therefore, sum of its residues at the points $Q_k^i$ equals zero, i.e.,
\beq\label{ny1}
\sum_{i,k} \tilde r_{k}^i\tilde\psi_0(t,Q_k^i)\tilde\psi_0^{\s}(t,Q_k^i)=0,
\ \ \tilde r_{k}^i=\res_{Q_{k}^i}d\Omega.
\eeq
Note, that equation (\ref{ny1}) is sufficient for the potential reduction
of the Schr\"odinger operator, which is equivalent to the equation $\tilde \psi_0(t,P_+)=1$.

From (\ref{dec1}) it follows that the evaluations of $\tilde\psi_0$ at the points
$Q_k^i$ satisfy linear equations
\beq\label{dec1a}
\sum_{i=1}^{\nu_k} \tilde c_{k,j}^{\,i}\tilde\psi_0(t,Q_k^i)=0, \ \ j=1,\ldots,n_k.
\eeq
Note, that the normalization of $\psi_0$ at $t=0$ implies
\beq\label{apr}
\sum_{i=1}^{\nu_k} \tilde c_{k,j}^{\,i}=0.
\eeq
The coefficients $\tilde c_{k,j}^i$ in (\ref{dec1a}) are
unique up to the transformations $\tilde c_{k,j}^i\to \sum_l g_{k,\,j}^{\,l}
\tilde c_{k,\,l}^i$, where $g_{\,k}=\{g_{k,\,j}^{\,l}\}$ are $t$-independent
non-degenerate matrices.
In what follows we normalize $\tilde c_k$ by the condition
\beq\label{normal}
\tilde c_{k,\,j}^{\,i}=\delta_j^i,\ i=1,\ldots,n_k.
\eeq
In that gauge the matrix elements $\tilde c_{k,j}^i(Z),\ i>n_k$ become well-defined
abelian functions on $X$.

Equations (\ref{ny1}) should follow from equations (\ref{dec1a}). The evaluations of
$\tilde\psi_0$ at the preimages of any two distinct singular points
$Q_k, Q_k', \ k\neq k'$ are independent. That implies
the following orthogonality relations
\beq\label{ny4}
 \sum_{j=1}^{n_k} \tilde r_{k}^j \tilde c_{k,j}^i \tilde c_{k,j}^{i}=-\tilde r_k^i,\ \ \
\sum_{j=1}^{n_k} \tilde r_{k}^j \tilde c_{k,j}^i \tilde c_{k,j}^{i'}=0, \ \
n_k<i\neq i'\leq\nu_k.
\eeq
(Compare (\ref{ny4}) with the orthogonality conditions established in (\cite{shiota2})).

\begin{cor} The multiplicity $\nu_k$ of the singular point $Q_k$ of the spectral curve is
equal to $\nu_k=2n_k$, where $n_k$ is the number of the linear relations in (\ref{dec1a}).
\end{cor}
{\it Proof.} As it was shown above, for the generic $Z$ the ring $\A^Z$ is maximal.
The ring $\A^Z$ is maximal if and only if for each $k$ the linear subspace
$W_k\subset \mathbb C^{2\nu_k}$ defined by the equations
$\sum_i \tilde c_{k,\,j}^{\,i}(Z)\,w_{k}^i=0$ is invariant under the
multiplication by a diagonal matrix $H_k$ {\it only} if $H_k$ is
a {\it scalar} matrix. The last condition
implies that each $(n_k\times n_k)$ minor of the matrix $c_k=\{\tilde c_{k,j}^{\,i}\}$
with $n_k<i\leq\nu_k$ is non-degenerate.
The columns of the matrix $\tilde c_{k,j}^{\,i}$ with $i>n_k$
are "orthogonal" to each other. Then, non-degeneracy of all the corresponding minors implies
that the number $\nu_k-n_k$ of such columns is not bigger than the dimension $n_k$
of the column vectors, i.e. $\nu_k\leq 2n_k$.

From (\ref{nnovp}) it follows that the degree of the pole divisor $D$ equals
$\deg D=\tilde g+1/2\sum_k \nu_k$, where $\tilde g$ is the genus of $\tilde \G$.
The uniqueness of the function $\psi_0$ defined by $D$ and the relations (\ref{dec1a}),
imply that $\deg D=\tilde g+\sum_k n_k$. The latter equations imply
$\sum_k(2n_k-\nu_k)=0$. As shown above, each term of the sum is non-negative.
Hence, $\nu_k=2n_k$ and the corollary
is thus proven.

\begin{lem} There exist constants $r_{k}^i$ such that the equation
\beq\label{ny3}
\sum_{i=1}^{2n_k} r_{k}^i \theta (A_k^i+Z)\, \theta (A_k^i- Z)=0, \ \ \ A_k^i=\tilde A(Q_k^i),
\eeq
holds.
\end{lem}
{\it Proof.}
Taking the square of (\ref{apr}) and using (\ref{ny4}) we get the equation
\beq\label{ny35}
\sum_{i=1}^{2n_k}\tilde r_k^i=0.
\eeq
The residues of the differential $d\Omega$ are well-defined  abelian functions
$\tilde r_{k}^i(Z)$ on $X$.
The pole divisor of all the residues coincides with the zero divisor of $\theta$.
The residue of $d\Omega$ at $Q_k^i$ equals zero, when the pole divisor of
$\tilde \psi_0$ contains the puncture $P_-$. Therefore, from (\ref{nnov51}) it follows that
$\tilde r_k^i$ has the form
\beq\label{ny2}
\tilde r_{k}^i(Z)=r_{k}^i{\theta (A_k^i+Z)\,\theta (A_k^i- Z)\over\theta^2(Z)}\,,
\eeq
where $r_{k}^i$ are constants. Equations (\ref{ny35}) and (\ref{ny2})
imply (\ref{ny3}) and the lemma is proved.

Our next and the final goal is to show that $n_k=1$, i.e. all of the singular points
of $\G$ are simple double points, as it is stated in the main
theorem.

If $n_k>1$, then from indecomposability of the matrix $\tilde c_{k,j}^{\,i}(Z)$ at the
generic $Z$ it follows that all the points $A_k^i$ are distinct, $A_k^i\neq A_k^j$.
That and the formula (\ref{nnov51}) for $\tilde \psi_0$ imply that
in the gauge (\ref{normal}) the coefficient $\tilde c_{k,j}^i$ for $i>n_k$ has pole at the
divisor $\theta (A_k^j+Z)=0$ and zero at the divisor $\theta(A_k^i+Z)=0$.

Let us fix a pair of indices $m,l>n_k$ and define a set $\D_k^{m,\,l}\subset X$
by the equations:
\beq\label{mar1}
\tilde c_{k,1}^i(Z)=0,\ \ n_k<i\neq m,l.
\eeq
On $\D_k^{m,\,l}$ equation (\ref{apr}) takes the form
\beq\label{apr1}
1+\tilde c_{k,1}^m(Z)+\tilde c_{k,1}^l(Z)=0,\ \ Z\in \D_k^{m,\,l}.
\eeq
From (\ref{mar1}) and the orthogonality conditions (\ref{ny4}) it follows  that
on $\D_k^{m,\,l}$ the equation $\tilde c_{k,1}^{\,l}(Z)=0$ implies
\beq\label{mar2}
\tilde c_{k,j}^{\,m}(Z)=0,\ \ j=n_k+2,\ldots,2n_k.
\eeq
Then, from (\ref{mar2}) it follows that
\beq\label{mar3}
Z\in \D_k^{m,\,l},\ \tilde c_{k,1}^{\,l}(Z)=0
\Rightarrow \tilde r_k^1+\tilde r_k^{m}=0.
\eeq
Hence, $\tilde c_{k,1}^{\,l}$, restricted to $\D_k^{m,\,l}$, is of the form
\beq\label{apr2}
\tilde c_{k,1}^{\,l}={\tilde r_k^1+\tilde r_k^m\over h(Z)}\, \cdot\,{\theta(A_k^{\,l}+Z)\over
\theta(A_k^1+Z)}\cdot\theta^2(Z)\,,\ \ Z\in \D_k^{m,\,l}\,,
\eeq
where $h$ is a {\it holomorphic} section of the line bundle
of $|2\Theta-A_k^{\,l}+A_k^1|$ restricted to $\D_k^{m,\,l}$. (Recall, that
$\tilde c_{k,j}^{\,i}$ and $\tilde r_k^i$ are abelian functions.)

The same arguments imply that on $\D_k^{m,\,l}$ zeros of $\tilde c_{k,1}^{\,m}$ are in the
zero divisor of $\tilde r_k^1+\tilde r_k^{\,l}$. Then using (\ref{apr1}) we get
\beq\label{apr3}
\tilde c_{k,1}^{\,m}=g{\tilde r_k^1+\tilde r_k^{\,l}\over h(Z)}\, \cdot\,{\theta(A_k^{\,m}+Z)\over
\theta(A_k^1+Z)}\cdot\theta^2(Z)\,,\ \ Z\in \D_k^{m,\,l}\,,
\eeq
where $g$ is a constant.
Therefore, $h$ is a section of the restriction to $\D_{k}^{m,\,l}$
of the line bundle $|2\Theta-A_k^{\,m}+A_k^1|$. Hence, $A_k^{m}=A_k^{\,l}$. The choice
of points $A_k^m,A_k^l$ was arbitrary, Therefore, we have
proved that all the points $A_k^i=A_k$ do coincide. In that case,
equations (\ref{dec1a}) are equivalent to $(2n_k-1)$ equations
$\tilde \psi_0(t,Q_k^i)=\tilde \psi_0(t,Q_k^j)$. That implies $2n_k-1=n_k=1$, i.e.,
$\G$ has only simple double singular points.
For such a curve all the sheafs $\F$ are line bundles. Therefore, the map
$j$ in (\ref{is}) is inverse to $i_Z$ in (\ref{imb}) and the main
theorem is thus proven.

\section{Appendix.}
\begin{theo} (R.Friedman) Let $\Gamma$ be an irreducible  projective curve, with an
involution $\sigma$, and suppose that the generalized Prym variety $P(\Gamma, \sigma)$
is compact. Then every singular point $x$ of $\Gamma$ is a fixed point of $\sigma$,
the singularity at $x$ is locally analytically isomorphic to a union
of coordinate axes in a neighborhood of the origin in $\Cee^N$ for some $N$,
and in a neighborhood of such a singular point $\sigma$ fixes
each of the local analytic branches.
\end{theo}
{\it Proof} Let $p\colon \widetilde{\Gamma}\to \Gamma$ be the normalization of $\Gamma$.
The involution $\sigma$ lifts to an involution on $\widetilde{\Gamma}$, also denoted
$\sigma$. The sheaf $p_*(\scrO_{\widetilde{\Gamma}})/\scrO_{\Gamma}$ is supported at
the finitely many points of $\Gamma_{\textrm{sing}}$.
The cohomology long exact sequence for
$$0\to \scrO_{\Gamma} \to p_*(\scrO_{\widetilde{\Gamma}})
\to p_*(\scrO_{\widetilde{\Gamma}})/\scrO_{\Gamma} \to 0$$ yields a long exact sequence
\begin{eqnarray}
0 \to H^0(\Gamma; \scrO_{\Gamma}) \to H^0(\Gamma; p_*(\scrO_{\widetilde{\Gamma}}))=
H^0(\widetilde{\Gamma};
\scrO_{\widetilde{\Gamma}})\to H^0(\Gamma ;
p_*(\scrO_{\widetilde{\Gamma}})/\scrO_{\Gamma}) \nonumber\\
\to H^1(\Gamma;\scrO_{\Gamma}) \to H^1(p_*(\scrO_{\widetilde{\Gamma}}))=
H^1(\widetilde{\Gamma}; \scrO_{\widetilde{\Gamma}}) \to 0. \nonumber
\end{eqnarray}
Since $\widetilde{\Gamma}$ is irreducible,
$H^0(\widetilde{\Gamma}; \scrO_{\widetilde{\Gamma}}) =\Cee$,
so that there is an exact sequence
$$0\to H^0(\Gamma ; p_*(\scrO_{\widetilde{\Gamma}})/\scrO_{\Gamma})
\to H^1(\Gamma;\scrO_{\Gamma}) \to  H^1(\widetilde{\Gamma}; \scrO_{\widetilde{\Gamma}})
\to 0. $$
Here $H^1(\Gamma;\scrO_{\Gamma})$ is the tangent space to the generalized
Jacobian of $\Gamma$ and the subspace
$H^0(\Gamma ; p_*(\scrO_{\widetilde{\Gamma}})/\scrO_{\Gamma})$
is the tangent space to its noncompact part.
If $V$ is a vector space on which $\sigma$ acts, let $V^{-}$
denotes the anti-invariant part of $V$, i.e.\ the $(-1)$-eigenspace.
Then the tangent space $T_{P(\Gamma, \sigma)}$ to $P(\Gamma, \sigma)$
fits into an exact sequence
$$0 \to H^0(\Gamma; p_*(\scrO_{\widetilde{\Gamma}})/\scrO_{\Gamma})^{-} \to T_{P(\Gamma, \sigma)} \to T_{P(\widetilde{\Gamma}, \sigma)}\to 0.$$
It follows that $P(\Gamma, \sigma)$ is compact if and only if
$H^0(\Gamma; p_*(\scrO_{\widetilde{\Gamma}})/\scrO_{\Gamma})^{-} = 0$.

First let us show that, for all $x\in \Gamma_{\textrm{sing}}$, $\sigma (x) = x$.
There is an isomorphism
$$H^0(\Gamma; p_*(\scrO_{\widetilde{\Gamma}})/\scrO_{\Gamma}) =
\bigoplus_{x\in \Gamma_{\rm sing}} \widetilde{R}_x/R_x,$$
where $R_x$ is the local ring $\scrO_{\Gamma, x}$
and $\widetilde{R}_x$ is its normalization.
Clearly, $\sigma$ induces  an isomorphism
$\widetilde{R}_x/R_x \cong \widetilde{R}_{\sigma(x)}/R_{\sigma(x)}$. If $\sigma(x)\neq x$, and $\alpha$ is a nonzero element of $\widetilde{R}_x/R_x$, then $\alpha - \sigma(\alpha) \in H^0(\Gamma; p_*(\scrO_{\widetilde{\Gamma}})/\scrO_{\Gamma})$ is nonzero, a contradiction. Hence $\sigma(x) = x$. Moreover, for all $\alpha \in \widetilde{R}_x/R_x$, $\sigma (\alpha) = \alpha$.

We now fix attention on a given $x\in \Gamma_{\textrm{sing}}$,
and write $R = R_x$ and $\widetilde{R} = \widetilde{R}_x$. Note that, if
$y_1, \dots, y_n$ are the preimages of $x$ in $\widetilde{\Gamma}$,
and $t_i$ is a local analytic coordinate for $\widetilde{\Gamma}$ at $y_i$, then
$\widetilde{R} \cong \bigoplus_i\Cee\{t_i\}$. Moreover, $R$ is a subalgebra of
$\widetilde{R}$, and $\dim_{\Cee} (\widetilde{R}/R) < \infty$; on particular,
$\widetilde{R}$ is a finite $R$-module. Let $\mathfrak{m}_i = t_i\Cee\{t_i\}$.
Clearly, $R$ is contained in the subalgebra $\Cee \oplus \bigoplus_i\mathfrak{m}_i$,
and the claim about the analytic nature of the singularities is just the statement
that $R = \Cee \oplus \bigoplus_i\mathfrak{m}_i$.

Next we claim that $\sigma$ does not permute the analytic branches through $x$.
If it did, then the action of $\sigma$ on $\widetilde{R}$ would exchange two
factors $\Cee\{t_i\}$ and $\Cee\{t_j\}$ for $j\neq i$. In this case, let $e_i$
be the image of $1\in \Cee\{t_i\} $ in $\widetilde{R}$. Then
$\sigma(e_i) - e_i\in (\widetilde{R})^{-}$, and $\sigma(e_i) - e_i \notin R$.
Thus $\sigma(e_i) - e_i$ is a nonzero class in
$(\widetilde{R}/R)^{-} = (\widetilde{R})^{-}/R^{-}$, a contradiction.
Thus, for every $i$, $\sigma$ fixes $\mathbb{C}\{t_i\}$ and induces  a holomorphic
involution on the corresponding branch of $\widetilde {\Gamma}$.

We can thus choose the coordinate $t_i$ so that $\sigma(t_i) = -t_i$.
Since $(\widetilde{R}/R)^{-} = (\widetilde{R})^{-}/R^{-} = 0$, $t_i\in R$ for every $i$.
Clearly, $t_i\in \mathfrak{m}$, where $\mathfrak{m}$ is the maximal ideal of the local
ring $R$. Now let
$\widehat{R} = \Cee \oplus \bigoplus_i\mathfrak{m}_i \subseteq \widetilde{R}$.
As noted above, $R\subseteq \widehat{R}$ and we must show that $R = \widehat{R}$.
In any case, $\widehat{R}$ is a finite $R$-module.  Given $r\in \widehat{R}$, there
exists a $c\in \Cee \subseteq R \subseteq \widehat{R}$ such that
$r-c\in \bigoplus_i\mathfrak{m}_i = (t_1, \dots, t_n)\widehat{R}\subseteq \mathfrak{m}
\widehat{R}$.
Thus $\widehat{R} = R + \mathfrak{m}\widehat{R}$, and so by Nakayama's
lemma $R = \widehat{R}$.

\medskip

{\footnotesize{\sl Acknowledgments. \rm The author wishes to thank Sam Grushevski for very
useful conversation on the subject of this paper, Robert Friedman for the communication
of the proof of Corollary 4.4. The author is grateful to Takahiro Shiota
whose remarks helped the author to clarify some missing arguments in
the first version of the paper.}}

\end{document}